\documentclass[a4paper,12pt]{article}

 \usepackage{amsfonts,amsmath,mathrsfs,amssymb}
 \usepackage{times,helvet,courier,type1cm}

 \usepackage{harvard}
 \usepackage{color}

 \allowdisplaybreaks

 \newtheorem{thm0}{Theorem}[section]
 \newtheorem{exa0}{Theorem}[section]

 \newtheorem{def1}[thm0]{Definition}
 \newtheorem{lem1}[thm0]{Lemma}
 \newtheorem{thm1}[thm0]{Theorem}
 \newtheorem{cor1}[thm0]{Corollary}
 \newtheorem{pro1}[thm0]{Proposition}
 \newtheorem{con1}[thm0]{Condition}
 \newtheorem{exa1}[exa0]{\it{Example}}

 \def\bglemma{\begin{lem1}}\def\edlemma{\end{lem1}}
 \def\bgtheorem{\begin{thm1}}\def\edtheorem{\end{thm1}}
 
 \def\bgproposition{\begin{pro1}}\def\edproposition{\end{pro1}}
 
 \def\bgexample{\begin{exa1}\rm{}\def\edexample{\end{exa1}}}

 \def\benumerate{\begin{enumerate}}\def\eenumerate{\end{enumerate}}
 \def\bitemize{\begin{itemize}}\def\eitemize{\end{itemize}}\def\itm{\item}

 \def\beqlb{\begin{eqnarray}}\def\eeqlb{\end{eqnarray}}
 \def\beqnn{\begin{eqnarray*}}\def\eeqnn{\end{eqnarray*}}

 \def\eqref#1{{\rm(\ref{#1})}}

 \def\ar{\!\!\!&}\def\nnm{\nonumber}\def\ccr{\nnm\\}

 \def\<{\langle}\def\>{\rangle}

 \def\mcr{\mathscr}\def\mbb{\mathbb}\def\mbf{\mathbf}
 \def\mrm{\mathrm}

 \def\proof{\noindent{\textit{Proof.~~}}}
 \def\qed{\hfill$\square$\smallskip}

 \def\d{\mrm{d}}\def\e{\mrm{e}}

 \def\qqquad{\qquad\qquad}

\begin{document}

\noindent{(Version: 2019-01-25)}

\bigskip\bigskip

\noindent{\Large\bf Sample paths of continuous-state branching}

\smallskip

\noindent{\Large\bf processes with dependent immigration}\,\footnote{Supported by the
National Natural Science Foundation of China (No.11531001).}

\bigskip

\noindent{Zenghu Li}

\medskip

\noindent{\it School of Mathematical Sciences, Beijing Normal University,}

\noindent{\it Beijing 100875, China. E-mail: lizh@bnu.edu.cn}

\bigskip

{\narrower

\noindent{\textit{Abstract:}} We prove the existence and pathwise uniqueness of the solution to a stochastic integral equation driven by Poisson random measures based on Kuznetsov measures for a continuous-state branching process. That gives a direct construction of the sample path of a continuous-state branching process with dependent immigration. The immigration rates depend on the population size via some functions satisfying a Yamada--Watanabe type condition. We only assume the existence of the first moment of the process. The existence of excursion law for the continuous-state branching process is not required. By special choices of the ingredients, we can make changes in the branching mechanism or construct models with competition.

\bigskip

\noindent{\textit{Key words:}} continuous-state branching process; dependent immigration; stochastic equation; Poisson random measure; Yamada--Watanabe type condition.

\smallskip

\noindent{\textit{MSC (2010) Subject Classification:}} 60J80, 60H10, 60H20

\par}

\bigskip

\section{Introduction}

 \setcounter{equation}{0}

The study of continuous-state branching processes (CB-processes) was started by Feller (1951), who noticed that a branching diffusion process may arise in a limit theorem of Galton--Watson discrete branching processes; see also Aliev and Shchurenkov (1982), Grimvall (1974) and Lamperti (1967a). A characterization of CB-processes by random time changes of L\'{e}vy processes was given by Lamperti (1967b). Continuous-state branching processes with immigration (CBI-processes) are natural generalizations of the CB-processes. The convergence of rescaled discrete branching processes with immigration to CBI-processes was studied in Aliev (1985), Kawazu and Watanabe (1971) and Li (2006, 2011). From a mathematical point of view, the continuous-state processes are usually easier to deal with because both their time and state spaces are smooth, and the distributions that appear are infinitely divisible. A continuous CBI-process with subcritical branching mechanism was used by Cox et al.\ (1985) to describe the evolution of interest rates and it has been known in mathematical finance as the \textit{Cox--Ingersoll--Ross model} (CIR-model). Compared with other financial models introduced before, the CIR-model is more appealing as it is positive  ($=$\,nonnegative) and mean-reverting. For general treatments and backgrounds of CB- and CBI-processes, the reader may refer to Kyprianou (2014) and Li (2011, 2018b). More complicated continuous-state population models involving a competition mechanism were studied in Pardoux (2016), which extend the stochastic logistic growth model of Lambert (2005).

In this paper, we are interested in a class of Markov processes, which we call \textit{continuous-state branching processes with dependent immigration} (CBDI-processes). By dependent immigration we mean the immigration rate depends on the state of the population via some function. This kind of immigration was studied in Dawson and Li (2003) and Fu and Li (2004) for measure-valued diffusions and extended in Li (2011) to general branching and immigration mechanisms. In those references, the immigration models were constructed in terms of stochastic equations driven by Poisson point measures on some path spaces. This approach is essential since the uniqueness of the corresponding martingale problems is usually unknown. In the references mentioned above, the immigration rate functions were assumed to be Lipschitz and the existence of some excursion laws of the corresponding measure-valued branching processes without immigration was required. In fact, some Poisson point measures based on the excursion laws were used there to represent the continuous part of the immigration.

The main purpose of this paper is to give a construction of the CBDI-process in terms of a stochastic equation of the type of Li (2011), but with non-Lipschitz immigration rate functions. A special case of the construction was given in the recent work of Li and Zhang (2019) by arguments of tightness and weak convergence, which require the existence of the second moment and the excursion law for the corresponding CB-process. We here replace the second moment assumption by the first moment one and remove the assumption on the existence of the excursion law. We focus on the one-dimensional model to simplify the presentations, but the arguments carry over to the measure-valued setting. In the stochastic equation considered here, the continuous immigration is represented by an increasing deterministic path and a Poisson point measure based on a Kuznetsov measure of the CB-process, which is slightly different from the equation in Li and Zhang (2019). The point of our approach is it gives a direct construction of the sample path of the CBDI-process with general branching and immigration mechanisms from those of the corresponding CB-process without immigration. By special choices of the ingredients, we can make changes in the branching mechanism of the CB-process or construct a model with competition. These kinds of constructions have been proved useful for the study of some financial problems; see, e.g., Bernis and Scotti (2018+) and Jiao et al.\ (2017). A more precise description of our results is given as follows.

Suppose that $c\ge 0$ and $b$ are real constants and $m(\d z)$ is a $\sigma$-finite measure on $(0,\infty)$ satisfying $\int_0^\infty (z\land z^2)m(\d z)< \infty$. Let $\phi$ be a function on $[0,\infty)$ defined by
 \beqlb\label{s1.1}
\phi(\lambda) = b\lambda + c\lambda^2 + \int_{(0,\infty)} (\e^{-z\lambda}-1+z\lambda) m(\d z).
 \eeqlb
A Markov process with state space $[0,\infty)$ is called a \textit{CB-process} with \textit{branching mechanism} $\phi$ if it has transition semigroup $(Q_t)_{t\ge 0}$ given by
 \beqlb\label{s1.2}
\int_{[0,\infty)} \e^{-\lambda y}Q_{t}(x,\d y)
 =
\e^{-x v_{t}(\lambda)}, \qquad \lambda\ge 0,~x\ge0,
 \eeqlb
where $t\mapsto v_{t}(\lambda)$ is the unique positive solution of
 \beqlb\label{s1.3}
\frac{\partial}{\partial t}v_{t}(\lambda) = -\phi(v_{t}(\lambda)),
 \qquad
v_{0}(\lambda)=\lambda.
 \eeqlb
The family of functions $(v_t)_{t\ge 0}$ satisfies $v_{r+t} = v_r\circ v_t$ for $r,t\ge 0$ and is called the \textit{cumulant semigroup} of the process. This semigroup has canonical representation
 \beqlb\label{s1.4}
v_t(\lambda) = h_t\lambda + \int_{(0,\infty)} (1-\e^{-\lambda z}) l_t(\d z),
\qquad t\ge 0,\lambda\ge 0,
 \eeqlb
where $h_t\ge 0$ and $l_t(\d z)$ is a $\sigma$-finite measure on $(0,\infty)$ satisfying $\int_{(0,\infty)} z l_t(\d z)< \infty$; see Li (2011, 2018b).

A generalization of the CB-process is described as follows. Let $\beta\ge 0$ be a constant and $\nu(\d z)$ a $\sigma$-finite measure on $(0,\infty)$ satisfying $\int_0^\infty (1\land z)\nu(\d z)< \infty$. Let $\psi$ be an \textit{immigration mechanism} defined by
 \beqlb\label{s1.5}
\psi(\lambda) = \beta\lambda + \int_{(0,\infty)} (1-\e^{-\lambda z}) \nu(\d z),
\qquad \lambda\ge 0.
 \eeqlb
A Markov process with state space $[0,\infty)$ is called a \textit{CBI-process} if it has transition semigroup $(P_t)_{t\ge 0}$ given by
 \beqlb\label{s1.6}
\int_{[0,\infty)} \e^{-\lambda y}P_{t}(x,\d y)
 =
\exp\bigg\{-x v_{t}(\lambda) - \int_0^t\psi(v_s(\lambda))\d s\bigg\},
 \eeqlb
where $\lambda\ge 0$ and $x\ge 0$.

From \eqref{s1.2} and \eqref{s1.6} we see that $(Q_t)_{t\ge 0}$ and $(P_t)_{t\ge 0}$ are Feller semigroups on $[0,\infty)$. Let $\mcr{D}= C^2[0,\infty)$ be the set of bounded continuous real functions on $[0,\infty)$ with bounded continuous derivatives up to the second order. By Theorem~9.30 in Li (2011), the generator $L_0$ of the CB-process and the generator $L_1$ of the CBI-process are respectively defined by
 \beqlb\label{s1.7}
L_0f(x) = cxf^{\prime\prime}(x) - bxf^\prime(x) + x\int_{(0,\infty)} \big[f(x+z) - f(x) - zf^\prime(x)\big] m(\d z),
 \eeqlb
and
 \beqlb\label{s1.8}
L_1f(x) = L_0f(x) + \beta f^{\prime}(x) + \int_{(0,\infty)} \big[f(x+z) - f(x)\big] \nu(\d z), \quad x\ge 0, f\in \mcr{D}.
 \eeqlb

We are interested in a generalization of the generators defined by \eqref{s1.7} and \eqref{s1.8}. Let $x\mapsto \beta(x)$ and $(x,z)\mapsto q(x,z)$ be positive Borel functions on $[0,\infty)$ and $[0,\infty)\times (0,\infty)$, respectively. We assume the following conditions:
 \bitemize

\itm[{\rm(1.A)}] (linear growth condition) there is a constant $K\ge 0$ so that
 \beqnn
\beta(x) + \int_{(0,\infty)} q(x,z)z\nu(\d z)\le K(1+x), \qquad x\ge 0;
 \eeqnn

\itm[{\rm(1.B)}] (Yamada--Watanabe type condition) there is an increasing and concave function $u\mapsto r(u)$ on $[0,\infty)$ so that $\int_{0+} r(u)^{-1} \d u= \infty$ and, for $x,y\ge 0$,
 \beqnn
|\beta(x)-\beta(y)| + \int_{(0,\infty)} |q(x,z)-q(y,z)|z \nu(\d z)\le r(|x-y|).
 \eeqnn

 \eitemize
By a \textit{CBDI-process} we mean a Markov process in $[0,\infty)$ with generator $L$ defined by, for $x\ge 0$ and $f\in \mcr{D}$,
 \beqlb\label{s1.9}
Lf(x) = L_0f(x) + \beta(x)f^\prime(x) + \int_{(0,\infty)} \big[f(x+z) - f(x)\big] q(x,z) \nu(\d z).
 \eeqlb

\bgexample\label{es1.1} When $\beta(x)\equiv \beta$ and $q(x,z)\equiv 1$ are constants, the operator $(L,\mcr{D})$ defined by \eqref{s1.9} generates a classical CBI-process. \edexample

\bgexample\label{es1.2} If $\beta(x)\equiv \beta x$ and $q(x,z)\equiv x$ for some constant $\beta\ge 0$, then the operator $(L,\mcr{D})$ defined by \eqref{s1.9} generates a CB-process with branching mechanism
 \beqnn
\lambda\mapsto \phi(\lambda) - \beta\lambda - \int_{(0,\infty)} (1-\e^{-\lambda z})\nu(\d z).
 \eeqnn
Then a change of the branching mechanism can be achieved by using dependent immigration. \edexample

\bgexample\label{es1.3} Let $x\mapsto G(x)$ be a positive function on $[0,\infty)$ satisfying the Yamada--Watanabe type condition and there is a constant $\beta> 0$ so that $G(x)\le \beta x$ for all $x\ge 0$. By setting $\beta(x)\equiv \beta x - G(x)$ and $q(x,z)\equiv 0$ in \eqref{s1.9} we get, for $x\ge 0$ and $f\in \mcr{D}$,
 \beqnn
Lf(x) = L_0f(x) + \beta xf^\prime(x) - G(x)f^\prime(x).
 \eeqnn
Then $(L,\mcr{D})$ generates a CB-process with competition; see, e.g., Berestycki et al.\ (2018), Lambert (2005) and Pardoux (2016). A more general class of population models, called \textit{continuous-state nonlinear branching processes}, have been studied in Li (2018a) and Li et al.\ (2017+). \edexample

Let $(\Omega, \mcr{F}, \mcr{F}_t, \mbf{P})$ be a filtered probability space satisfying the usual conditions. Let $\{B(t)\}$ be a $(\mcr{F}_t)$-Brownian motion and let $\{M(\d s,\d z,\d u)\}$ and $\{N(\d s,\d z,\d u)\}$ be $(\mcr{F}_t)$-Poisson random measures on $(0,\infty)^3$ with intensities $\d sm(\d z)\d u$ and $\d s\nu(\d z)\d u$, respectively. Suppose that $\{B(t)\}$, $\{M(\d s,\d z,\d u)\}$ and $\{N(\d s,\d z,\d u)\}$ are independent of each other. Let $\{\tilde{M}(\d s,\d z,\d u)\}$ be the compensated measure of $\{M(\d s,\d z,\d u)\}$. By Theorem~5.1 in Fu and Li (2010), for any $\mcr{F}_0$-measurable positive random variable $y(0)$ there is a pathwise unique positive solution to the stochastic equation:
 \beqlb\label{s1.10}
y(t) \ar=\ar y(0) + \int_0^t \sqrt{2cy(s-)} \d B(s) + \int_0^t\int_{0}^{\infty} \int_{0}^{y(s-)} z \tilde{M}(\d s,\d z,\d u)\cr
 \ar\ar
+ \int_0^t [\beta(y(s-))-by(s-)] \d s +\int_0^t\int_{0}^{\infty}\int_{0}^{q(y(s-),z)} z N(\d s,\d z,\d u).
 \eeqlb
Here and in the sequel, we understand that, for $b\ge a\ge 0$,
 \beqlb\label{s1.11}
\int_a^b = -\int_b^a = \int_{(a,b]} ~\mbox{and}~ \int_a^\infty = \int_{(a,\infty)}.
 \eeqlb
The reader may refer to Ikeda and Watanabe (1989) and Situ (2005) for the theory of stochastic equations. We will show that a positive c\`{a}dl\`{a}g process $\{y(t): t\ge 0\}$ is a weak solution of \eqref{s1.10} if and only if it solves the martingale problem of $(L,\mcr{D})$. Then \eqref{s1.10} gives a construction of the CBDI-process. From this equation we see that the immigration of $\{y(t): t\ge 0\}$ involves two parts: the \textit{continuous part} given by the drift $\beta(y(s-))\d s$ and the \textit{discontinuous part} determined by the intensity $q(y(s-),z)$ and Poisson random measure $N(\d s,\d z,\d u)$. Stochastic equations in forms similar to \eqref{s1.10} have also been studied in Bertoin and Le~Gall (2006), Dawson and Li (2006, 2012) and Li and Ma (2015).

Let $D[0,\infty)$ denote the space of positive c\`{a}dl\`{a}g paths $w = \{w(t): t\ge 0\}$. For any $w\in D[0,\infty)$ let $\alpha(w)= \inf\{s\ge 0: w(s)> 0\}$ and $\zeta(w)= \sup\{s\ge 0: w(s)>0\}$. Let $W$ be the set of paths $w\in D[0,\infty)$ such that $w(t)>0$ for $\alpha(w)< t< \zeta(w)$ and $w(t)=0$ for $t<\alpha(w)$ or $t\ge \zeta(w)$. Let $[0]\in W$ be the path that is constantly zero. On the space $W$ we define the $\sigma$-algebras $\mcr{W}= \sigma(\{w(s): s\ge 0\})$ and $\mcr{W}_t = \sigma(\{w(s): 0\le s\le t\})$ for $t\ge 0$.

From \eqref{s1.2} we see that zero is a trap for the CB-process. Let $(Q_t^\circ)_{t\ge 0}$ be the restriction of the transition semigroup $(Q_t)_{t\ge 0}$ on $(0,\infty)$. For a $\sigma$-finite measure $\mu$ on $(0,\infty)$ write
 \beqnn
\mu Q_t^\circ(\d y) = \int_{(0,\infty)} \mu(\d x) Q_t^\circ(x,\d y), \qquad
t\ge 0, y> 0.
 \eeqnn
A family of $\sigma$-finite measures $(\kappa_t)_{t>0}$ on $(0,\infty)$ is called an \textit{entrance rule} for $(Q^\circ_t)_{t\ge 0}$ if $\kappa_r Q^\circ_{t-r}\le \kappa_t$ for all $t>r>0$ and $\kappa_rQ^\circ_{t-r}\to \kappa_t$ as $r\to t$.

Let $(l_t)_{t>0}$ be the family of $\sigma$-finite measures on $(0,\infty)$ determined by \eqref{s1.4}. By Theorems~3.13 and~3.15 in Li (2018b) we see that $(l_t)_{t>0}$ is an entrance rule for $(Q_t^\circ)_{t\ge 0}$, which is referred to as the \textit{canonical entrance rule} of the CB-process. We shall give a simple and direct construction of the $\sigma$-finite measure $\mbf{N}_0$ on $(W,\mcr{W})$ such that, for $0<t_1< t_2< \cdots< t_n$ and $x_1,x_2,\ldots,x_n\in (0,\infty)$,
 \beqlb\label{s1.12}
\ar\ar\mbf{N}_0(\alpha(w)\le t_1,w(t_1)\in \d x_1, w(t_2)\in \d x_2, \ldots, w(t_n)\in \d x_n,t_n< \zeta(w)) \ccr
 \ar\ar\qqquad
=\, l_{t_1}(\d x_1) Q_{t_2-t_1}^\circ(x_1,\d x_2) Q_{t_3-t_2}^\circ(x_2,\d x_3) \cdots Q_{t_n-t_{n-1}}^\circ(x_{n-1},\d x_n). \quad
 \eeqlb
We call $\mbf{N}_0$ the \textit{canonical Kuznetsov measure} of the CB-process. The existence of this measure is also a consequence of a general result on Markov processes; see Theorem~3.8 in Glover and Getoor (1987, p.63). The expression \eqref{s1.12} intuitively means that the coordinate process $\{w(t): t>0\}$ under $\mbf{N}_0$ is a Markov process in $(0,\infty)$ with transition semigroup $(Q_t^\circ)_{t\ge 0}$ and one-dimensional distributions $(l_t)_{t>0}$.

Let $\mbf{Q}_x$ denote the distribution of the CB-process $\{x(t): t\ge 0\}$ with initial value $x(0)=x\ge 0$. Let $\{(X_t,\mcr{F}_t)\}$ be a CB-process with generator $L_0$ defined by \eqref{s1.7} with $\mbf{P}[X_0]< \infty$. Let $\{N_0(\d s,\d u,\d w)\}$ be a Poisson random measure on $(0,\infty)^2\times W$ with intensity $\d s\d u\mbf{N}_0(\d w)$ and $\{N_1(\d s,\d z,\d u,\d w)\}$ a Poisson random measure on $(0,\infty)^3\times W$ with intensity $\d s\nu(\d z)$ $\d u\mbf{Q}_z(\d w)$. Suppose that $\{(X_t,\mcr{F}_t)\}$, $\{N_0(\d s,\d u,\d w)\}$ and $\{N_1(\d s,\d z,\d u,\d w)\}$ are defined on a complete probability space and are independent of each other. For $t\ge 0$ let $\mcr{G}_t = \sigma(\mcr{F}_t\cup \mcr{H}_t)$, where
 \beqnn
\mcr{H}_t \ar=\ar \sigma(\{N_0((0,s]\times B\times A), N_1((0,s]\times C\times A): \cr
 \ar\ar\qqquad
0<s\le t, B\in \mcr{B}(0,\infty), C\in \mcr{B}(0,\infty)^2, A\in \mcr{W}_{t-s}\}).
 \eeqnn
We consider the stochastic integral equation
 \beqlb\label{s1.13}
Y_t \ar=\ar X_t + \int_0^t h_{t-s}\beta(Y_{s-})\d s + \int_0^t\int_0^{\beta(Y_{s-})} \int_W w(t-s) N_0(\d s,\d u,\d w) \cr
 \ar\ar\qquad
+ \int_0^t\int_0^\infty\int_0^{q(Y_{s-},z)} \int_W w(t-s) N_1(\d s,\d z,\d u,\d w).
 \eeqlb
By a solution of \eqref{s1.13} we mean a positive c\`{a}dl\`{a}g process $\{Y_t: t\ge 0\}$ that is adapted to the filtration $(\mcr{G}_t)$ and satisfies the equation a.s.\ for each $t\ge 0$.

The main result of this paper shows that there is a pathwise unique solution $\{Y_t: t\ge 0\}$ of \eqref{s1.13} and $\{(Y_t,\mcr{G}_t): t\ge 0\}$ is indeed a CBDI-process. Here the second and third terms on the right-hand side of the equation represent the continuous part of the immigration and the last term represents the discontinuous immigration. In Section~2, we provide the direct construction of the canonical Kuznetsov measure defined by \eqref{s1.12}. We also give a reformulation of the Markov property of the measure that is more convenient for our applications. In Section~3, we construct some inhomogeneous immigration processes with deterministic immigration rates. In Section~4, we construct immigration processes with predictable immigration rates and prove some useful properties of them. The existence and uniqueness of solution to the stochastic equation \eqref{s1.13} are established in Section~5.

\section{The canonical Kuznetsov measure}

 \setcounter{equation}{0}

In this section, we give a simple and direct construction of the canonical Kuznetsov measure defined by \eqref{s1.12}. For notational convenience, we extend the definition of each path $w\in W$ by setting $w(s)=0$ for $s<0$. Recall that $\mbf{Q}_x$ denotes the distribution on $(W,\mcr{W})$ of the CB-process $\{x(t): t\ge 0\}$ with $x(0)=x\ge 0$. For $r\ge 0$ and $w\in W$, we define $\rho_rw\in W$ by $\rho_rw(t) = w(t-r)$. Let $\mbf{Q}_x^{(r)}$ denote the image of $\mbf{Q}_x$ induced by the map $w\mapsto \rho_rw$. Then $\mbf{Q}_x^{(r)}$ is supported by $\{w\in W: \alpha(w)=r, w(r)= x\}$. Given any $\sigma$-finite measure $\mu$ on $[0,\infty)$ let
 \beqlb\label{s2.1}
\mbf{Q}_\mu^{(r)}(A) = \int_{[0,\infty)} \mbf{Q}_x^{(r)}(A)\mu(\d x),
 \qquad
A\in \mcr{W}.
 \eeqlb
In particular, if $\mu$ is a probability measure, then $\mbf{Q}_\mu^{(0)}$ is the distribution on $(W,\mcr{W})$ of the CB-process with initial distribution $\mu$ and $\mbf{Q}_\mu^{(r)}$ is the image of $\mbf{Q}_\mu^{(0)}$ induced by the map $w\mapsto \rho_rw$.

\bglemma\label{ts2.1} For any $s> r\ge 0$ and any positive $\mcr{W}_r$-measurable function $F$ on $W$ we have $\mbf{Q}_\mu^{(s)}(\{w\in W: F(w)\neq F([0])\})= 0$. \edlemma

\proof Let $0\le r_1\le r_2\le\dots\le r_n\le r< s$ and let $f_1,f_2,\dots,f_n$ be bounded Borel functions on $[0,\infty)$. For any $x\ge 0$ we have
 \beqnn
\ar\ar\mbf{Q}_x^{(s)}\big[f_1(w(r_1))f_2(w(r_2))\cdots f_n(w(r_n))\big] \ccr
 \ar\ar\qquad
= \mbf{Q}_x\big[f_1(w(r_1-s))f_2(w(r_2-s))\cdots f_n(w(r_n-s))\big] \ccr
 \ar\ar\qquad
= f_1(0)f_2(0)\cdots f_n(0) \ccr
 \ar\ar\qquad
= f_1([0](r_1))f_2([0](r_2))\cdots f_n([0](r_n)).
 \eeqnn
A monotone class argument shows $\mbf{Q}_x^{(s)}F(w)= F([0])$ for any positive $\mcr{W}_r$-measurable function $F$ on $W$, and so $\mbf{Q}_x^{(s)}f(F(w))= f(F([0]))$ for any bounded Borel function $f$ on $[0,\infty)$. Then we have $\mbf{Q}_x^{(s)}(\{w\in W: F(w)\neq F([0])\})= 0$, which implies the desired result by (\ref{s2.1}). \qed

\bgtheorem\label{ts2.2} Let $(l_t)_{t>0}$ be the canonical entrance rule for the CB-process determined by \eqref{s1.4}. Then there is a unique $\sigma$-finite measure $\mbf{N}_0$ on $(W,\mcr{W})$ that does not charge the singleton $\{[0]\}\in \mcr{W}$ and satisfies \eqref{s1.12}. Moreover, we have:
 \benumerate

\itm[(1)] If $\phi^\prime(\infty)= \infty$, then $\mbf{N}_0$ is supported by $\{w\in W: \alpha(w)=0, w(0)=0\}$.

\itm[(2)] If $\delta:= \phi^\prime(\infty)< \infty$, then $\mbf{N}_0$ is supported by $\{w\in W: \alpha(w)>0, w(\alpha(w))>0\}$ and has the representation
 \beqlb\label{s2.2}
\mbf{N}_0(F) = \int_0^\infty \e^{-\delta s} \mbf{Q}_m^{(s)}(F)\d s, \qquad F\in \mcr{W}.
 \eeqlb
 \eenumerate
\edtheorem

\proof In the case $\phi^\prime(\infty)= \infty$, we have $h_t=0$ for all $t>0$ and the result follows by Theorem~6.1 in Li (2018b). In the case $\delta:= \phi^\prime(\infty)< \infty$, let us define the $\sigma$-finite measure $\mbf{N}_0$ on $(W,\mcr{W})$ by \eqref{s2.2}. Since $\mbf{Q}_m^{(r)}$ is supported by $\{w\in W: \alpha(w)=r, w(r)>0\}$, we see that $\mbf{N}_0$ is supported by $\{w\in W: \alpha(w)>0, w(\alpha(w))>0\}$. Let $f_1,\ldots,f_n$ be positive Borel functions on $[0,\infty)$ with $f_1(0) =\ldots =f_n(0) =0$. Then
 \beqnn
\mbf{N}_0[f_1(w(t_1))]
 \ar=\ar
\mbf{N}_0\big[\alpha(w)\le t_1, f_1(w(t_1))\big] \ccr
 \ar=\ar
\int_0^{t_1} \e^{-\delta s} \mbf{Q}_m^{(s)}[f_1(w(t_1))]\d s \cr
 \ar=\ar
\int_0^{t_1} \e^{-\delta s} \mbf{Q}_m^{(0)}[f_1(w(t_1-s))]\d s \cr
 \ar=\ar
\int_0^{t_1}\e^{-\delta s}\d s\int_{(0,\infty)}f_1(y) mQ_{t_1-s}^\circ(\d y) \cr
 \ar=\ar
\int_{(0,\infty)}f_1(y) l_{t_1}(\d y),
 \eeqnn
where the last equality follows by Theorem~3.15 in Li (2018b). Then \eqref{s1.12} holds for $n=1$. From the Markov property of $\mbf{Q}_m^{(s)}$ it follows that, for $n\ge 2$,
 \beqnn
\ar~\ar\mbf{N}_0\Big[\alpha(w)\le t_1, f_1(w(t_1))\cdots f_{n-1}(w(t_{n-1}))f_n(w(t_n)), t_n< \zeta(w)\Big] \cr
 \ar~\ar\qquad
= \mbf{N}_0\Big[\alpha(w)\le t_1, f_1(w(t_1))\cdots f_{n-1}(w(t_{n-1}))f_n(w(t_n))\Big] \cr
 \ar~\ar\qquad
= \int_0^{t_1} \e^{-\delta s} \mbf{Q}_m^{(s)}\Big[f_1(w(t_1))\cdots f_{n-1}(w(t_{n-1})) f_n(w(t_n))\Big]\d s \cr
 \ar~\ar\qquad
= \int_0^{t_1} \e^{-\delta s} \mbf{Q}_m^{(0)}\Big[f_1(w(t_1-s))\cdots f_{n-1}(w(t_{n-1}-s)) f_n(w(t_n-s))\Big]\d s \cr
 \ar~\ar\qquad
= \int_0^{t_1} \e^{-\delta s}\mbf{Q}_m^{(0)}\Big[f_1(w(t_1-s))\cdots f_{n-1}(w(t_{n-1}-s)) Q_{t_n-t_{n-1}}^\circ f_n(w(t_{n-1}-s))\Big]\d s \cr
 \ar~\ar\qquad
= \int_0^{t_1} \e^{-\delta s}\mbf{Q}_m^{(s)}\Big[f_1(w(t_1))\cdots f_{n-1}(w(t_{n-1})) Q_{t_n-t_{n-1}}^\circ f_n(w(t_{n-1}))\Big]\d s \cr
 \ar~\ar\qquad
= \mbf{N}_0\Big[f_1(w(t_1))\cdots f_{n-1}(w(t_{n-1}))Q_{t_n-t_{n-1}}^\circ f_n(w(t_{n-1}))\Big].
 \eeqnn
Then we get \eqref{s1.12} by induction, which determines the measure $\mbf{N}_0$ uniquely by the measure extension theorem. \qed

\bgtheorem\label{ts2.3} Let $t\ge r> 0$ and let $F$ be a positive $\mcr{W}_r$-measurable function on $W$. Then for any $\lambda\ge 0$ we have
 \beqlb\label{s2.3}
\mbf{N}_0[F(w)(1-\e^{-\lambda w(t)})]
 =
\mbf{N}_0[F(w)(1-\e^{-v_{t-r}(\lambda)w(r)})] + F([0])[h_rv_{t-r}(\lambda) - h_t\lambda].
 \eeqlb
\edtheorem

\proof In the case $\phi^\prime(\infty)= \infty$, we have $h_r=h_t=0$ and \eqref{s2.3} follows from \eqref{s1.12}. In the case $\delta:= \phi^\prime(\infty)< \infty$, the measure $\mbf{N}_0$ is given by \eqref{s2.2}. By Lemma~\ref{ts2.1}, for any $s>r$ we have $\mbf{Q}_m^{(s)}(\{w\in W: F(w)\neq F([0])\})= 0$. Then we use the Markov property of $\mbf{Q}_m^{(s)}$ to see
 \beqnn
\mbf{N}_0[F(w)(1-\e^{-\lambda w(t)})]
 \ar=\ar
\mbf{N}_0\big[\alpha(w)\le t, F(w)(1-\e^{-\lambda w(t)})\big] \ccr
 \ar=\ar
\int_0^t \e^{-\delta s}\mbf{Q}_m^{(s)}[F(w)(1-\e^{-\lambda w(t)})]\d s \cr
 \ar=\ar
\int_0^r \e^{-\delta s}\mbf{Q}_m^{(s)}[F(w)(1-\e^{-v_{t-r}(\lambda)w(r)})]\d s \cr
 \ar\ar\qquad\qquad
+ \int_r^t \e^{-\delta s}\mbf{Q}_m^{(s)}[F([0])(1-\e^{-\lambda w(t)})]\d s \cr
 \ar=\ar
\mbf{N}_0\big[F(w)(1-\e^{-v_{t-s}(\lambda)w(r)})\big] \ccr
 \ar\ar\qquad\qquad
+\, F([0])\int_r^t \e^{-\delta s}\d s\int_{(0,\infty)} (1-\e^{-xv_{t-s}(\lambda)})m(\d x) \ccr
 \ar=\ar
\mbf{N}_0[F(w)(1-\e^{-v_{t-s}(\lambda)w(r)})] + F([0])\e^{-\delta r}[v_{t-r}(\lambda)-\e^{-\delta(t-r)}\lambda],
 \eeqnn
where the last equality holds by Theorem~3.15 of Li (2018b). \qed

The existence of Markovian measures determined by entrance rules was first noticed by Kuznetsov (1974). In the setting of Borel right Markov processes, it was proved in Glover and Getoor (1987). The relation \eqref{s2.3} gives a reformulation of the Markov property of the canonical Kuznetsov measure, which is more convenient than \eqref{s1.12} for the application in the next section. In the special case of $\phi^\prime(\infty)= \infty$, we have $h_t=0$ for every $t>0$ and the canonical Kuznetsov measure $\mbf{N}_0$ is known as the \textit{excursion law} of the CB-process. In that special case, the property \eqref{s2.3} is already known; see, e.g., the proof of Theorem~8.24 in Li (2011).

\section{Deterministic immigration rates}

 \setcounter{equation}{0}

In this section, we give constructions of some inhomogeneous CBI-processes with deterministic immigration rates from random paths selected by Poisson point measures. The reader may refer to Li (2011, 2018b) for similar constructions. Let $\{(X_t,\mcr{F}_t)\}$, $\{N_0(\d s,\d u,\d w)\}$ and $\{N_1(\d s,\d z,\d u,\d w)\}$ be as in the introduction.

\bgtheorem\label{ts3.1} Let $s\mapsto \rho(s)$ be a positive locally integrable function on $[0,\infty)$. For $t\ge 0$ let
 \beqlb\label{s3.1}
Y_t = X_t + \int_0^t h_{t-s}\rho(s)\d s + \int_0^t\int_0^{\rho(s)}\int_{W} w(t-s) N_0(\d s,\d u,\d w)
 \eeqlb
and let $\mcr{G}_t = \sigma(\mcr{F}_t\cup \mcr{H}_t^0)$, where
 \beqlb\label{s3.2}
\mcr{H}_t^0 = \sigma(\{N_0((0,s]\times B\times A): 0<s\le t, B\in \mcr{B}(0,\infty), A\in \mcr{W}_{t-s}\}).
 \eeqlb
Then $\{(Y_t,\mcr{G}_t): t\ge 0\}$ is a Markov process in $[0,\infty)$ with inhomogeneous transition semigroup $(P_{r,t}^\rho)_{t\ge r\ge 0}$ given by
 \beqlb\label{s3.3}
\int_{[0,\infty)}\e^{-\lambda y}P_{r,t}^\rho(x,\d y)
 =
\exp\bigg\{-xv_{t-r}(\lambda) - \int_r^tv_{t-s}(\lambda)\rho(s)\d s\bigg\}.
 \eeqlb
\edtheorem

\proof Let $Z_t = Y_t - X_t$ for $t\ge 0$. We first prove that $\{(Z_t,\mcr{H}^0_t)\}$ is a Markov process with transition semigroup $(P_{r,t}^\rho)_{t\ge r\ge 0}$. Let $t\ge r\ge \tau\ge 0$ and let $F$ be a positive function on $(0,\infty)^2\times W$ measurable with respect to $\mcr{B}(0,\tau]\times \mcr{B}(0,\infty)\times \mcr{W}_{r-\tau}$. It suffices to show, for any $\lambda\ge 0$,
 \beqnn
\ar\ar\mbf{P}\bigg[\exp\bigg\{-\int_0^\tau\int_0^\infty\int_{W} F(s,u,w) N_0(\d s,\d u,\d w) - \lambda Z_t\bigg\}\bigg] \cr
 \ar\ar\qquad
=\, \mbf{P}\bigg[\exp\bigg\{-\int_0^\tau\int_0^\infty\int_{W} F(s,u,w) N_0(\d s,\d u,\d w) \cr
 \ar\ar\qqquad\qqquad\qqquad
-\, v_{t-r}(\lambda)Z_r - \int_r^tv_{t-s}(\lambda)\rho(s)\d s\bigg\}\bigg].
 \eeqnn
Since $\{1_{\{w=[0]\}} N_0(\d s,\d u,\d w)\}$ is independent of $\{1_{\{w\neq [0]\}} N_0(\d s,\d u,\d w)\}$ and $\{Z_t\}$, we can assume $F(s,u,[0])=0$ in the following calculations. Writing $G(s,u,w)= F(s,u,w)1_{\{s\le \tau\}} + \lambda w(t-s)1_{\{u\le \rho(s)\}}$, we have
 \beqnn
\ar\ar\mbf{P}\bigg[\exp\bigg\{-\int_0^\tau\int_0^\infty\int_{W} F(s,u,w) N_0(\d s,\d u,\d w) - \lambda Z_t\bigg\}\bigg] \cr
 \ar\ar\qquad
=\, \mbf{P}\bigg[\exp\bigg\{-\lambda\int_0^t h_{t-s}\rho(s)\d s - \int_0^t\int_0^\infty \int_{W} G(s,u,w) N_0(\d s,\d u,\d w)\bigg\}\bigg] \cr
 \ar\ar\qquad
=\, \exp\bigg\{-\lambda\int_0^t h_{t-s}\rho(s)\d s - \int_0^t\d s\int_0^\infty \mbf{N}_0(1-\e^{-G(s,u,w)}) \d u\bigg\} \cr
 \ar\ar\qquad
=\, \exp\bigg\{-\lambda\int_0^t h_{t-s}\rho(s)\d s - \int_0^t\d s\int_0^\infty \mbf{N}_0(1-\e^{-F(s,u,w) 1_{\{s\le \tau\}}}) \d u\bigg\} \cr
 \ar\ar\qquad\quad
\cdot\,\exp\bigg\{-\int_0^t\d s\int_0^\infty \mbf{N}_0\big[\e^{-F(s,u,w)1_{\{s\le \tau\}}}(1-\e^{-\lambda w(t-s)1_{\{u\le \rho(s)\}}})\big] \d u\bigg\} \cr
 \ar\ar\qquad
=\, \exp\bigg\{-\lambda\int_0^t h_{t-s}\rho(s)\d s - \int_0^t\d s\int_0^\infty \mbf{N}_0(1-\e^{-F(s,u,w) 1_{\{s\le \tau\}}}) \d u\bigg\} \cr
 \ar\ar\qquad\quad
\cdot\,\exp\bigg\{-\int_0^r\d s\int_0^{\rho(s)} \mbf{N}_0\big[\e^{-F(s,u,w) 1_{\{s\le \tau\}}}(1-\e^{-\lambda w(t-s)})\big] \d u\bigg\} \cr
 \ar\ar\qquad\quad
\cdot\,\exp\bigg\{-\int_r^t\d s\int_0^{\rho(s)} \mbf{N}_0(1-\e^{-\lambda w(t-s)}) \d u\bigg\},
 \eeqnn
where, by Theorem~\ref{ts2.3},
 \beqnn
\ar\ar \mbf{N}_0\big[\e^{-F(s,u,w) 1_{\{s\le \tau\}}}(1-\e^{-\lambda w(t-s)})\big] \ccr
 \ar\ar\qquad
= \mbf{N}_0\big[\e^{-F(s,u,w)1_{\{s\le \tau\}}} (1-\e^{-v_{t-r}(\lambda)w(r-s)})\big] + h_{r-s}v_{t-r}(\lambda)-h_{t-s}\lambda.
 \eeqnn
Then we can use \eqref{s1.4} and continue the calculation
 \beqnn
\ar\ar\mbf{P}\bigg[\exp\bigg\{-\int_0^\tau\int_0^\infty\int_{W} F(s,u,w) N_0(\d s,\d u,\d w) - \lambda Z_t\bigg\}\bigg] \cr
 \ar\ar\qquad
=\, \exp\bigg\{-\lambda\int_0^r h_{t-s}\rho(s)\d s - \lambda\int_r^t h_{t-s}\rho(s)\d s\bigg\} \cr
 \ar\ar\qquad\quad
\cdot\,\exp\bigg\{- \int_0^r\d s\int_0^\infty \mbf{N}_0(1- \e^{-F(s,u,w)1_{\{s\le \tau\}}})\d u\bigg\} \cr
 \ar\ar\qquad\quad
\cdot\,\exp\bigg\{-\int_0^r\d s\int_0^{\rho(s)} \mbf{N}_0\big[\e^{-F(s,u,w)1_{\{s\le \tau\}}}(1-\e^{-v_{t-r}(\lambda)w(r-s)})\big]\d u\bigg\} \cr
 \ar\ar\qquad\quad
\cdot\,\exp\bigg\{-\int_0^r\big[h_{r-s}v_{t-r}(\lambda)-h_{t-s}\lambda\big]\rho(s) \d s\bigg\} \cr
 \ar\ar\qquad\quad
\cdot\,\exp\bigg\{- \int_r^t\rho(s)\d s\int_{(0,\infty)} (1-\e^{-y\lambda}) l_{t-s}(\d y)\bigg\} \cr
 \ar\ar\qquad
=\, \exp\bigg\{-v_{t-r}(\lambda)\int_0^r h_{r-s}\rho(s)\d s - \int_r^tv_{t-s}(\lambda)\rho(s)\d s\bigg\} \cr
 \ar\ar\qquad\quad
\cdot\,\exp\bigg\{- \int_0^r\d s\int_0^\infty \mbf{N}_0(1- \e^{-F(s,u,w)1_{\{s\le \tau\}}})\d u\bigg\} \cr
 \ar\ar\qquad\quad
\cdot\,\exp\bigg\{-\int_0^r\d s\int_0^\infty \mbf{N}_0\big[\e^{-F(s,u,w)1_{\{s\le \tau\}}}(1-\e^{-v_{t-r}(\lambda)w(r-s)1_{\{u\le \rho(s)\}}})\big]\d u\bigg\} \cr
 \ar\ar\qquad
=\, \exp\bigg\{-v_{t-r}(\lambda)\int_0^r h_{r-s}\rho(s)\d s - \int_r^tv_{t-s}(\lambda)\rho(s)\d s\bigg\} \cr
 \ar\ar\qquad\quad
\cdot\,\exp\bigg\{-\int_0^r\d s\int_0^\infty\mbf{N}_0\big(1-\e^{-F(s,u,w) 1_{\{s\le \tau\}} - v_{t-r}(\lambda)w(r-s)1_{\{u\le \rho(s)\}}}\big)\d u\bigg\} \cr
 \ar\ar\qquad
=\, \exp\bigg\{-v_{t-r}(\lambda)\int_0^r h_{r-s}\rho(s)\d s - \int_r^tv_{t-s}(\lambda)\rho(s)\d s\bigg\} \cr
 \ar\ar\qquad\quad
\cdot\,\mbf{P}\bigg[\exp\bigg\{-\int_0^r\int_0^\infty\int_{W} \big[F(s,u,w) 1_{\{s\le \tau\}} \cr
 \ar\ar\qqquad\qqquad\qqquad
+\, v_{t-r}(\lambda) w(r-s)1_{\{u\le \rho(s)\}}\big] N_0(\d s,\d u,\d w)\bigg\}\bigg] \cr
 \ar\ar\qquad
=\, \exp\bigg\{-v_{t-r}(\lambda)\int_0^r h_{r-s}\rho(s)\d s - \int_r^tv_{t-s}(\lambda)\rho(s)\d s\bigg\} \cr
 \ar\ar\qquad\quad
\cdot\,\mbf{P}\bigg[\exp\bigg\{-\int_0^\tau\int_0^\infty\int_{W} F(s,u,w) N_0(\d s,\d u,\d w) \cr
 \ar\ar\qqquad\qqquad\qqquad
-\int_0^r\int_0^{\rho(s)}\int_{W} v_{t-r}(\lambda) w(r-s) N_0(\d s,\d u,\d w)\bigg\}\bigg] \cr
 \ar\ar\qquad
=\, \mbf{P}\bigg[\exp\bigg\{-\int_0^\tau\int_0^\infty\int_{W} F(s,u,w) N_0(\d s,\d u,\d w) \cr
 \ar\ar\qqquad\qqquad\qqquad
-\, v_{t-r}(\lambda)Z_r - \int_r^tv_{t-s}(\lambda)\rho(s)\d s\bigg\}\bigg].
 \eeqnn
Then $\{(Z_t,\mcr{H}^0_t)\}$ is a Markov process with transition semigroup $(P_{r,t}^\rho)_{t\ge r\ge 0}$. Using the independence of $\{(X_t,\mcr{F}_t)\}$ and $\{(Z_t,\mcr{H}^0_t)\}$, one can see $\{(Y_t,\mcr{G}_t)\}$ is a Markov process with transition semigroup $(P_{r,t}^\rho)_{t\ge r\ge 0}$. \qed

\bgtheorem\label{ts3.2} Let $s\mapsto \rho(s)$ be a positive locally integrable function on $[0,\infty)$ and $(s,z)\mapsto g(s,z)$ a positive measurable function on $[0,\infty)\times (0,\infty)$ such that
 \beqnn
\int_0^t\d s\int_{(0,\infty)} g(s,z)z \nu(\d z)< \infty, \qquad t\ge 0.
 \eeqnn
For $t\ge 0$ let
 \beqlb\label{s3.4}
Y_t \ar=\ar X_t + \int_0^t h_{t-s}\rho(s)\d s + \int_0^t\int_0^{\rho(s)}\int_{W} w(t-s) N_0(\d s,\d u,\d w) \cr
 \ar~\ar\qquad
+ \int_0^t\int_0^\infty\int_0^{g(s,z)}\int_{W} w(t-s) N_1(\d s,\d z,\d u,\d w)
 \eeqlb
and let $\mcr{G}_t = \sigma(\mcr{F}_t\cup \mcr{H}_t^0\cup \mcr{H}_t^1)$, where $\mcr{H}_t^0$ is defined by \eqref{s3.2} and
 \beqnn
\mcr{H}_t^1 = \sigma(\{N_1((0,s]\times C\times A): 0<s\le t, C\in \mcr{B}(0,\infty)^2, A\in \mcr{W}_{t-s}\}).
 \eeqnn
Then $\{(Y_t,\mcr{G}_t): t\ge 0\}$ is a Markov process in $[0,\infty)$ with inhomogeneous transition semigroup $(P_{r,t}^{\rho,g})_{t\ge r\ge 0}$ given by
 \beqlb\label{s3.5}
\int_{[0,\infty)}\e^{-\lambda y}P_{r,t}^{\rho,g}(x,\d y)
 \ar=\ar
\exp\bigg\{-xv_{t-r}(\lambda) - \int_r^t v_{t-s}(\lambda)\rho(s)\d s \cr
 \ar~\ar\qquad
- \int_r^t\d s\int_{(0,\infty)} (1-\e^{-zv_{t-s}(\lambda)})g(s,z)\nu(\d z)\bigg\}.
 \eeqlb
\edtheorem

\proof Let $Z_t$ denote the last term on the right-hand side of \eqref{s3.4}. Let $t\ge r\ge \tau\ge 0$ and let $F$ be a positive function on $(0,\infty)^3\times W$ measurable relative to $\mcr{B}[0,\tau]\times \mcr{B}(0,\infty)^2\times \mcr{W}_{r-\tau}$. For $\lambda\ge 0$, writing $H(s,z,u,w)= F(s,z,u,w) 1_{\{s\le \tau\}} + \lambda w(t-s) 1_{\{u\le g(s,z)\}}$, we have
 \beqnn
\ar\ar\mbf{P}\bigg[\exp\bigg\{-\int_0^\tau\int_0^\infty\int_0^\infty\int_{W} F(s,z,u,w)N_1(\d s,\d z,\d u,\d w) - \lambda Z_t\bigg\}\bigg] \cr
 \ar\ar\quad
=\,\mbf{P}\bigg[\exp\bigg\{-\int_0^t\int_0^\infty\int_0^\infty\int_{W} H(s,z,u,w) N_1(\d s,\d z,\d u,\d w)\bigg\}\bigg] \cr
 \ar\ar\quad
=\,\exp\bigg\{-\int_0^t\d s\int_0^\infty\nu(\d z)\int_0^\infty \mbf{Q}_z(1-\e^{-H(s,z,u,w)}) \d u\bigg\} \cr
 \ar\ar\quad
=\,\exp\bigg\{-\int_0^t\d s\int_0^\infty\nu(\d z)\int_0^\infty\mbf{Q}_z(1-\e^{-F(s,z,u,w) 1_{\{s\le \tau\}}})\d u \cr
 \ar\ar\qqquad
-\int_0^t\d s\int_0^\infty\nu(\d z)\int_0^{g(s,z)} \mbf{Q}_z\big[\e^{-F(s,z,u,w) 1_{\{s\le \tau\}}}(1-\e^{-\lambda w(t-s)})\big]\d u\bigg\} \cr
 \ar\ar\quad
=\,\exp\bigg\{-\int_0^r\d s\int_0^\infty\nu(\d z)\int_0^\infty\mbf{Q}_z(1-\e^{-F(s,z,u,w) 1_{\{s\le \tau\}}})\d u \cr
 \ar\ar\qqquad
-\int_0^r\d s\int_0^\infty\nu(\d z)\int_0^{g(s,z)} \mbf{Q}_z\big[\e^{-F(s,z,u,w) 1_{\{s\le \tau\}}}(1-\e^{-v_{t-r}(\lambda)w(r-s)})\big]\d u \cr
 \ar\ar\qqquad
-\int_r^t\d s\int_0^\infty\nu(\d z)\int_0^{g(s,z)}\mbf{Q}_z(1-\e^{-\lambda w(t-s)})\d u\bigg\} \cr
 \ar\ar\quad
=\,\mbf{P}\bigg[\exp\bigg\{-\int_0^r\int_0^\infty\int_0^\infty\int_{W} \big[F(s,z,u,w)1_{\{s\le \tau\}} \ccr
 \ar\ar\qqquad\qqquad
+\, v_{t-r}(\lambda)w(t-s)1_{\{u\le g(s,z)\}}\big]N_1(\d s,\d z,\d u,\d w) \ccr
 \ar\ar\qqquad
-\int_r^t\d s\int_0^\infty \big(1-\e^{-zv_{t-s}(\lambda)}\big)g(s,z)\nu(\d z)\bigg\}\bigg] \cr
 \ar\ar\quad
=\,\mbf{P}\bigg[\exp\bigg\{-\int_0^\tau\int_0^\infty\int_0^\infty\int_{W} F(s,z,u,w) N_1(\d s,\d z,\d u,\d w) \cr
 \ar\ar\qqquad
-\, v_{t-r}(\lambda)Z_r - \int_r^t\d s\int_0^\infty \big(1-\e^{-zv_{t-s}(\lambda)}\big)g(s,z)\nu(\d z)\bigg\}\bigg].
 \eeqnn
Then $\{(Z_t,\mcr{H}_t^1)\}$ is a Markov process in $[0,\infty)$ with inhomogeneous transition semigroup $(P_{r,t}^{0,g})_{t\ge r\ge 0}$ defined by \eqref{s3.5} with $\rho=0$. By Theorem~\ref{ts3.1} we see $\{(Y_t-Z_t,\sigma(\mcr{F}_t\cup \mcr{H}_t^0))\}$ is a Markov process in $[0,\infty)$ with inhomogeneous transition semigroup $(P_{r,t}^{\rho,0})_{t\ge r\ge 0}$ defined by \eqref{s3.5} with $g=0$. Then the desired result follows by the independence of those two processes. \qed

We may think of the process $\{Y_t: t\ge 0\}$ defined by \eqref{s3.4} as an inhomogeneous CBI-process with \textit{immigration rates} given by $\{(\rho(s),g(s,z)): s\ge 0,z>0\}$.

\section{Predictable immigration rates}

 \setcounter{equation}{0}

Suppose that $\{(X_t,\mcr{F}_t)\}$, $\{N_0(\d s,\d u,\d w)\}$ and $\{N_1(\d s,\d z,\d u,\d w)\}$ are given as in the introduction. Let the filtration $(\mcr{G}_t)$ be defined as in Theorem~\ref{ts3.2}. Let $\mcr{L}^1$ denote the set of $(\mcr{G}_t)$-predictable processes $\rho = \{\rho(t): t\ge 0\}$ satisfying
 \beqnn
\|\rho\|_t := \mbf{P}\bigg[\int_0^t|\rho(s)| \d s\bigg]< \infty,
\qquad t\ge 0.
 \eeqnn
We identify $\rho_1, \rho_2\in \mcr{L}^1$ if $\|\rho_1-\rho_2\|_t = 0$ for every $t\ge 0$ and define the metric $d$ on $\mcr{L}^1$ by
 \beqnn
d(\rho_1,\rho_2) = \sum_{n=1}^\infty \frac{1}{2^n}(1\land \|\rho_1-\rho_2\|_n).
 \eeqnn
Let $\mcr{L}^1_\nu(0,\infty)$ denote the set of two-parameter processes $g = \{g(t,z): t\ge 0, z>0\}$ that are $(\mcr{G}_t)$-predictable in the sense of Li (2011, p.163) and satisfy
 \beqnn
\|g\|_{\nu,t} := \mbf{P}\bigg[\int_0^t\int_{(0,\infty)}|g(s,z)| z\nu(\d z)\d s\bigg]< \infty,
\qquad t\ge 0.
 \eeqnn
We identify $g_1, g_2\in \mcr{L}^1_\nu(0,\infty)$ if $\|g_1-g_2\|_{\nu,t} = 0$ for every $t\ge 0$ and define the metric $d_\nu$ on $\mcr{L}^1_\nu(0,\infty)$ by
 \beqnn
d_\nu(g_1,g_2) = \sum_{n=1}^\infty \frac{1}{2^n}(1\land \|g_1-g_2\|_{\nu,n}).
 \eeqnn

For $\rho\in \mcr{L}^1$, $g\in \mcr{L}^1_\nu(0,\infty)$ and a positive c\`{a}dl\`{a}g process $\{Y_t: t\ge 0\}$ we consider the following properties:
 \bitemize

\itm[{\rm(4.A)}] The process $\{Y_t: t\ge 0\}$ has no negative jumps and the optional random measure
 \beqnn
N_0(\d s,\d z) := \sum_{s>0}1_{\{\Delta Y_s\ne 0\}}\delta_{(s,\Delta Y_s)}(\d s,\d z),
 \eeqnn
where $\Delta Y_s = Y_s-Y_{s-}$, has predictable compensator
 \beqnn
\hat{N}_0(\d s,\d z)= Y_{s-}\d sm(\d z) + g(s,z)\d s\nu(\d z)
 =
Y_s\d sm(\d z) + g(s,z)\d s\nu(\d z).
 \eeqnn
Let $\tilde{N}_0(\d s,\d z) = N_0(\d s,\d z) - \hat{N}_0(\d s,\d z)$. We have
 \beqlb\label{s4.1}
Y_t = Y_0 + M^c(t) + M^d(t) + \int_0^t \bigg[\rho(s) + \int_{(0,\infty)} g(s,z)z\nu(\d z) - bY_s\bigg]\d s,
 \eeqlb
where $\{M^c(t): t\ge 0\}$ is a square-integrable continuous martingale with quadratic variation $2cY_{s-}\d s= 2cY_s\d s$ and
 \beqlb\label{s4.2}
M^d(t) = \int_0^t\int_0^\infty z\tilde{N}_0(\d s,\d z), \qquad t\ge 0,
 \eeqlb
is a purely discontinuous martingale.

\itm[{\rm(4.B)}] For every $f\in \mcr{D}$, we have
 \beqlb\label{s4.3}
f(Y_t)
 \ar=\ar
f(Y_0) + \int_0^t \big[\rho(s)f^\prime(Y_s) - bY_sf^\prime(Y_s) + c Y_s f^{\prime\prime}(Y_s)\big]\d s
\cr
 \ar\ar
+ \int_0^t\d s \int_{(0,\infty)} \big[f(Y_s+z) - f(Y_s)\big]g(s,z)\nu(\d z)
+ \mbox{local mart.}
\cr
 \ar\ar
+ \int_0^t Y_s\d s \int_{(0,\infty)} \big[f(Y_s+z) - f(Y_s) - zf^\prime(Y_s)\big]m(\d z).
 \eeqlb

\eitemize

\bgproposition\label{ts4.1} The above properties {\rm(4.A)} and {\rm(4.B)} are equivalent. \edproposition

\proof If $\{Y_t\}$ has property (4.A), we may use It\^{o}'s formula to see it has property (4.B). Conversely, let us assume $\{Y_t\}$ has property (4.B). For any $\lambda>0$, by applying this property to the function $f(x)=\e^{-\lambda x}$ we have
 \beqlb\label{s4.4}
\e^{-\lambda Y_t}
 \ar=\ar
\e^{-\lambda Y_0} + \int_0^t \e^{-\lambda Y_s}\big[bY_s\lambda - \rho(s)\lambda + c Y_s \lambda^2\big]\d s
\cr
 \ar\ar
+ \int_0^t g(s,z)\d s \int_{(0,\infty)} \e^{-\lambda Y_s}\big(\e^{-\lambda z} - 1\big)\nu(\d z) + \mbox{local mart.}
\cr
 \ar\ar
+ \int_0^t Y_s\d s \int_{(0,\infty)} \e^{-\lambda Y_s}\big(\e^{-\lambda z} - 1 + \lambda z\big)m(\d z).
 \eeqlb
Then the strictly positive process $\{\e^{-\lambda Y_t}: t\ge 0\}$ is a special semi-martingale. By It\^o's formula, we see $\{Y_t: t\ge 0\}$ is also a special semi-martingale. Now define an optional random measure $N_0(\d s,\d z)$ on $(0,\infty)\times \mbb{R}$ by
 \beqnn
N_0(\d s,\d z) = \sum_{s>0}1_{\{\Delta Y_s\ne 0\}}\delta_{(s,\Delta Y_s)}(\d s,\d z),
 \eeqnn
where $\Delta Y_s = Y_s-Y_{s-}$. Let $\hat{N}_0(\d s,\d z)$ denote the predictable compensator of $N_0(\d s,\d z)$ and let $\tilde{N}_0(\d s,\d z)$ denote the compensated random measure; see Dellacherie and Meyer (1982, pp.375). We can write
 \beqlb\label{s4.5}
Y_t = Y_0 + U_t + M^c_t + M^d_t,
 \eeqlb
where $\{U(t)\}$ is a predictable process with locally bounded variations, $\{M^c_t\}$ is a continuous local martingale and
 \beqnn
M^d(t) = \int_0^t\int_{\mbb{R}} z \tilde{N}_0(\d s,\d z)
 \eeqnn
is a purely discontinuous local martingale; see Dellacherie and Meyer (1982, p.353 and p.376) or Jacod and Shiryaev (2003, pp.84--85). Let $\{C_t\}$ denote the quadratic variation process of $\{M^c_t\}$. By \eqref{s4.5} and It\^{o}'s formula,
 \beqlb\label{s4.6}
\e^{-\lambda Y_t} \ar=\ar \e^{-\lambda Y_0} - \lambda\int_0^t\e^{-\lambda Y_{s-}}\d U(s)
+ \frac{1}{2}\lambda^2\int_0^t \e^{-\lambda Y_{s-}}\d C_s \cr
 \ar~\ar
+ \int_0^t\int_\mbb{R} \e^{-\lambda Y_{s-}} \big(\e^{-z\lambda} - 1 +
z\lambda\big) N_0(\d s,\d z) + \mbox{local mart.} \cr
 \ar=\ar
\e^{-\lambda Y_0} - \lambda\int_0^t\e^{-\lambda Y_{s-}}\d U(s)
+ \frac{1}{2}\lambda^2\int_0^t \e^{-\lambda Y_{s-}}\d C_s \cr
 \ar~\ar
+ \int_0^t\int_\mbb{R} \e^{-\lambda Y_{s-}} \big(\e^{-z\lambda} - 1 +
z\lambda\big) \hat{N}_0(\d s,\d z) + \mbox{local mart.}
 \eeqlb
But, the canonical decomposition of the special semi-martingale $\{\e^{-\lambda Y_t}: t\ge 0\}$ is unique; see, e.g., Dellacherie and Meyer (1982, p.213). From \eqref{s4.4} and \eqref{s4.6} we see $\d C_s = 2cY_s\d s$,
 \beqnn
\hat{N}_0(\d s,\d z)= Y_s\d sm(\d z) + g(s,z)\d s\nu(\d z)
 \eeqnn
and
 \beqnn
\d U(s)= \bigg[\rho(s) - bY_s + \int_{(0,\infty)} g(s,z)z \nu(\d z)\bigg]\d s.
 \eeqnn
Then the process $\{Y_t\}$ has no negative jumps. \qed

Now let us consider a kind of stochastic immigration rates. Given positive processes $\rho\in \mcr{L}^1$ and $g\in \mcr{L}^1_\nu(0,\infty)$, we define
 \beqlb\label{s4.7}
Y_t \ar=\ar X_t + \int_0^t h_{t-s}\rho(s)\d s + \int_0^t\int_0^{\rho(s)}\int_{W} w(t-s) N_0(\d s,\d u,\d w) \cr
 \ar\ar\qquad
+ \int_0^t\int_0^\infty\int_0^{g(s,z)}\int_{W} w(t-s) N_1(\d s,\d z,\d u,\d w).
 \eeqlb

\bgproposition\label{ts4.2} Let $\{Y_t: t\ge 0\}$ be defined by \eqref{s4.7}. Then for any $t\ge 0$ we have
 \beqlb\label{s4.8}
\mbf{P}[Y_t] = \e^{-bt}\mbf{P}[X_0] + \int_0^t \e^{-b(t-s)}\mbf{P}\bigg[\rho(s) + \int_0^\infty g(s,z)z \nu(\d z)\bigg]\d s.
 \eeqlb
\edproposition

\proof Recall that both $s\mapsto \rho(s)$ and $(s,z)\mapsto g(s,z)$ are predictable. From \eqref{s4.7} we have
 \beqnn
Y_t \ar=\ar X_t + \int_0^t h_{t-s}\rho(s)\d s + \int_0^t\int_0^{\rho(s)} \int_{W} w(t-s) \tilde{N}_0(\d s,\d u,\d w) \cr
 \ar\ar\qquad
+ \int_0^t\rho(s)\d s\int_{W} w(t-s) \mbf{N}_0(\d w) \cr
 \ar\ar\qquad
+ \int_0^t\int_0^\infty\int_0^{g(s,z)}\int_{W} w(t-s) \tilde{N}_1(\d s,\d z,\d u,\d w) \cr
 \ar\ar\qquad
+ \int_0^t\d s\int_0^\infty g(s,z)\mbf{Q}_z[w(t-s)]\nu(\d z).
 \eeqnn
By (3.7) in Li (2011) or (3.5) in Li (2018b) we have $\mbf{P}[X_t] = \mbf{P}[X_0]\e^{-bt} = \mbf{P}[Y_0]\e^{-bt}$. It follows that
 \beqnn
\mbf{P}[Y_t] \ar=\ar \mbf{P}[X_t] + \mbf{P}\bigg[\int_0^t h_{t-s}\rho(s)\d s\bigg] + \mbf{P}\bigg[\int_0^t\mbf{N}_0[w(t-s)]\rho(s)\d s\bigg] \cr
 \ar\ar\qquad
+\, \mbf{P}\bigg[\int_0^t\d s\int_0^\infty\mbf{Q}_z[w(t-s)]g(s,z)\nu(\d z)\bigg] \cr
 \ar=\ar
\e^{-bt}\mbf{P}[X_0] + \mbf{P}\bigg[\int_0^t\bigg(h_{t-s} + \int_0^\infty y l_{t-s}(\d y)\bigg)\rho(s)\d s\bigg] \cr
 \ar\ar\qquad
+\, \mbf{P}\bigg[\int_0^t \e^{-b(t-s)}\d s\int_0^\infty g(s,z)z\nu(\d z)\bigg].
 \eeqnn
Then we get desired equality from (3.4) in Li (2018b). \qed

We may interpret the process $\{Y_t: t\ge 0\}$ defined by \eqref{s4.7} as a generalization of the inhomogeneous CBI-process with \textit{predictable immigration rates} given by $\{(\rho(s),g(s,z)): s\ge 0, z> 0\}$. Recall that we identify $\rho_1, \rho_2\in \mcr{L}^1$ if $\|\rho_1-\rho_2\|_t = 0$ for every $t\ge 0$ and identify $g_1, g_2\in \mcr{L}^1_\nu(0,\infty)$ if $\|g_1-g_2\|_{\nu,t} = 0$ for every $t\ge 0$. By Proposition~\ref{ts4.2} one can see that choosing different representatives of $\rho\in \mcr{L}^1$ and $g\in \mcr{L}^1_\nu(0,\infty)$ in \eqref{s4.7} only gives different modifications of the process $\{Y_t: t\ge 0\}$.

\bgtheorem\label{ts4.3} The process $\{Y_t: t\ge 0\}$ defined by \eqref{s4.7} has a c\`{a}dl\`{a}g modification and it satisfies properties {\rm(4.A)} and {\rm(4.B)}. \edtheorem

The proof of the above theorem is based on approximations of $\rho\in \mcr{L}^1$ and $g\in \mcr{L}^1_\nu(0,\infty)$ using simpler processes. Let $\mcr{L}^0$ denote the set of processes $\rho\in \mcr{L}^1$ of the form
 \beqlb\label{s4.9}
\rho(s) = \rho(r_1)1_{\{0\}}(s) + \sum_{i=0}^\infty \rho(r_{i+1})1_{(r_i,r_{i+1}]}(s),
 \eeqlb
where $\{0= r_0< r_1< r_2< \cdots\}$ is a sequence increasing to infinity and each $\omega\mapsto \rho(\omega,r_{i+1})$ is $\mcr{G}_{r_i}$-measurable. Let $\mcr{L}^0_\nu(0,\infty)$ denote the set of processes $g\in \mcr{L}^1_\nu(0,\infty)$ of the form
 \beqlb\label{s4.10}
g(s,z) = g(r_1,z)1_{\{0\}}(s) + \sum_{i=0}^\infty g(r_{i+1},z)1_{(r_i,r_{i+1}]}(s),
 \eeqlb
where $\{0= r_0< r_1< r_2< \cdots\}$ is as above and each $(\omega,z)\mapsto g(\omega,r_{i+1},z)$ is $\mcr{G}_{r_i}\times \mcr{B}(0,\infty)$-measurable.

\bglemma\label{ts4.4} For positive processes $\rho\in \mcr{L}^0$ and $g\in \mcr{L}^0_\nu(0,\infty)$, the results of Theorem~\ref{ts4.3} hold. \edlemma

\proof By choosing $\{0= r_0< r_1< r_2< \cdots\}$ suitably, we can represent $s\mapsto \rho(s)$ and $(s,z)\mapsto g(s,z)$ by \eqref{s4.9} and \eqref{s4.10} using the same sequence. Then under $\mbf{P}(\cdot | \mcr{G}_{r_{i-1}})$ we can think of $\rho(r_i)$ as a deterministic constant and $z\mapsto g(r_i,z)$ as a deterministic function. For $i\ge 1$ let
 \beqnn
Y_i(t) \ar=\ar X_t + \int_0^{t\land r_i} h_{t-s}\rho(s)\d s + \int_0^{t\land r_i}\int_0^{\rho(s)} \int_W w(t-s) N_0(\d s,\d u,\d w) \cr
 \ar\ar\qqquad
+ \int_0^{t\land r_i}\int_0^\infty\int_0^{g(s,z)} \int_W w(t-s) N_1(\d s,\d z,\d u,\d w)
 \eeqnn
and $\mcr{G}_t^i = \sigma(\mcr{F}_t\cup \mcr{D}_t^i)$, where
 \beqnn
\mcr{D}_t^i \ar=\ar \sigma(\{N_0([(0,s\land r_i]\times B\times A), N_1((0,s\land r_i]\times C\times A): \cr
 \ar\ar\qqquad
0< s\le t, B\in \mcr{B}(0,\infty), C\in \mcr{B}(0,\infty)^2, A\in \mcr{W}_{t-s}\}).
 \eeqnn
Then $\{Y_i(t): t\ge 0\}$ is adapted to the filtration $\{\mcr{G}_t^i: t\ge 0\}$. Note also that $Y_i(t) = Y_t$, $\mcr{G}_t^i = \mcr{G}_t$ for $0\le t\le r_i$ and $Y_i(t)\le Y_t$, $\mcr{G}_t^i\subset \mcr{G}_t$ for $t\ge r_i$. We claim that the following properties hold:
 \bitemize

\itm[(a.1)] $\{(Y_i(t), \mcr{G}_t): r_{i-1}\le t\le r_i\}$ is a CBI-process under $\mbf{P}(\cdot |\mcr{G}_{r_{i-1}})$ with time-independent immigration rate $\{(\rho(r_i),g(r_i,z)): z>0\}$;

\itm[(a.2)] $\{(Y_i(t), \mcr{G}_t): t\ge r_i\}$ and $\{(Y_i(t), \mcr{G}_t^i): t\ge r_i\}$ are CB-processes under both $\mbf{P}(\cdot |\mcr{G}_{r_{i-1}})$ and $\mbf{P}(\cdot |\mcr{G}_{r_i})$.

 \eitemize
For $i=1$, those properties follow immediately from Theorem~\ref{ts3.2}. Suppose they hold for some $i\ge 1$. Let
 \beqnn
Z_i(t) \ar=\ar \int_{t\land r_i}^{t\land r_{i+1}} h_{t-s}\rho(s)\d s + \int_{t\land r_i}^{t\land r_{i+1}} \int_0^{\rho(s)} \int_W w(t-s) N_0(\d s,\d u,\d w) \cr
 \ar\ar\qquad
+ \int_{t\land r_i}^{t\land r_{i+1}}\int_0^\infty \int_0^{g(s,z)} \int_W w(t-s) N_1(\d s,\d z,\d u,\d w) \cr
 \eeqnn
and
 \beqnn
\mcr{H}_t^i \ar=\ar \sigma(\{N_0((r_i,s\land r_{i+1}]\times B\times A), N_1((r_i,s\land r_{i+1}]\times F\times A): \cr
 \ar\ar\qqquad
r_i<s\le t, B\in \mcr{B}(0,\infty), F\in \mcr{B}(0,\infty)^2, A\in \mcr{W}_{t-s}\}).
 \eeqnn
Using Theorem~\ref{ts2.3} again we see:

 \bitemize

\itm[(b.1)] $\{(Z_i(t), \mcr{H}_t^i): r_i\le t\le r_{i+1}\}$ is a CBI-process under $\mbf{P}(\cdot |\mcr{G}_{r_i})$ with time-independent immigration rate $\{\rho(r_{i+1}), g(r_{i+1},z)): z>0\}$;

\itm[(b.2)] $\{(Z_i(t), \mcr{H}_t^i): t\ge r_{i+1}\}$ is a CB-process under both $\mbf{P}(\cdot|\mcr{G}_{r_i})$ and $\mbf{P}(\cdot|\mcr{G}_{r_{i+1}})$.

 \eitemize
Here the processes $\{(Y_i(t), \mcr{G}_t^i): t\ge r_i\}$ and $\{(Z_i(t), \mcr{H}_t^i): t\ge r_i\}$ are independent of each other under $\mbf{P}(\cdot |\mcr{G}_{r_i})$. Note also that $Y_{i+1}(t) = Y_i(t) + Z_i(t)$ and $\mcr{G}_t^{i+1} = \sigma(\mcr{G}_t^i\cup \mcr{H}_t^i)$ for $t\ge r_i$. By Proposition~5.7 in Li (2018b), properties (a.1) and (a.2) also hold when $i$ is replaced by $i+1$. Then they hold for all $i\ge 1$ by induction. By applying Theorem~7.2 in Li (2018b) step by step on the intervals $[r_{i-1},r_i]$, $i=1,2,\cdots$ we see $\{Y_t: t\ge 0\}$ satisfies property {\rm(4.B)}. By Proposition~\ref{ts4.1} it also satisfies property {\rm(4.A)}. \qed

\bglemma\label{ts4.5} Suppose that $\rho\in \mcr{L}^1$ and $g\in \mcr{L}^1_\nu(0,\infty)$ are positive processes. Let $\{\rho_k\}\subset \mcr{L}^0$ and $\{g_k\}\subset \mcr{L}^0_\nu(0,\infty)$ be positive sequences such that $d(\rho_k,\rho) + d_\nu(g_k,g)\to 0$ as $k\to \infty$. Let $\{Y_k(t): t\ge 0\}$ be the positive c\`{a}dl\`{a}g process defined by (\ref{s4.7}) with $\rho=\rho_k$ and $g=g_k$. Then there is a positive c\`{a}dl\`{a}g process $\{Y(t): t\ge 0\}$ so that
 \beqlb\label{s4.11}
\lim_{k\to\infty} \mbf{P}\Big[\sup_{0\le s\le t} |Y_k(s) - Y(s)|\Big]= 0, \qquad t\ge 0,
 \eeqlb
and there is a subsequence $\{k_n\}\subset \{k\}$ so that a.s.\
 \beqlb\label{s4.12}
\lim_{n\to \infty} \sup_{0\le s\le t}|Y_{k_n}(s) - Y(s)| = 0, \qquad t\ge 0.
 \eeqlb
\edlemma

\proof For any $j,k\ge 1$, we can represent $\rho_j,\rho_k$ and $g_j,g_k$ in the form of \eqref{s4.9} and \eqref{s4.10} using the same sequence $\{0= r_0< r_1< r_2< \cdots\}$. Then we have $|Y_j(t) - Y_k(t)|\le Z_{j,k}(t)$, where
 \beqlb\label{s4.13}
Z_{j,k}(t) \ar=\ar \int_0^t h_{t-s}|\rho_j(s)-\rho_k(s)|\d s + \int_0^t \int_{\rho_j(s)\land \rho_k(s)}^{\rho_j(s)\vee \rho_k(s)} \int_W w(t-s) N_0(\d s,\d u,\d w) \cr
 \ar\ar\qquad
+ \int_0^t\int_0^\infty\int_{g_j(s,z)\land g_k(s,z)}^{g_j(s,z)\vee g_k(s,z)} \int_W w(t-s) N_1(\d s,\d z,\d u,\d w).
 \eeqlb
We can rewrite the above expression into
 \beqnn
Z_{j,k}(t) \ar=\ar \int_0^t h_{t-s}|\rho_j(s)-\rho_k(s)|\d s + \int_0^t \int_0^{|\rho_j(s)-\rho_k(s)|} \int_W w(t-s) N_0^{j,k}(\d s,\d u,\d w) \cr
 \ar\ar\qquad
+ \int_0^t\int_0^\infty\int_0^{|g_j(s,z)-g_k(s,z)|} \int_W w(t-s) N_1^{j,k}(\d s,\d z,\d u,\d w),
 \eeqnn
where $N_0^{j,k}(\d s,\d u,\d w)= N_0(\d s,\rho_j(s)\land \rho_k(s)+\d u,\d w)$ is a Poisson random measure with intensity $\d s\d u\mbf{N}_0(\d w)$ and $N_1^{j,k}(\d s,\d z,\d u,\d w)= N_1(\d s,\d z,g_j(s,z)\land g_k(s,z)+\d u,\d w)$ is a Poisson random measure with intensity $\d s\nu(\d z)\d u\mbf{Q}_z(\d w)$. One can see that $\{(Z_{j,k}(t), \mcr{G}_t): t\ge 0\}$ is a CBI-process with predictable immigration rates given by $\{(|\rho_j(s) - \rho_k(s)|, |g_j(s,z) - g_k(s,z)|): s\ge 0, z>0\}$. By Proposition~\ref{ts4.2} we see that
 \beqlb\label{s4.14}
\mbf{P}[Z_{j,k}(t)]\le \e^{|b|t}(\|\rho_j-\rho_k\|_t+\|g_j-g_k\|_{\nu,t}).
 \eeqlb
By Lemma~\ref{ts4.4}, the results of Theorem~\ref{ts4.3} hold for $\{(Z_{j,k}(t), \mcr{G}_t)\}$. By Proposition~\ref{ts4.1}, this process has properties (4.A) and (4.B). In particular, it has no negative jumps and the optional random measure
 \beqnn
N_{j,k}(\d s,\d z) = \sum_{s>0}1_{\{\Delta Z_{j,k}(s)\ne 0\}}\delta_{(s,\Delta Z_{j,k}(s))}(\d s,\d z),
 \eeqnn
where $\Delta Z_{j,k}(s)= Z_{j,k}(s) - Z_{j,k}(s-)$, has predictable compensator
 \beqnn
\hat{N}_{j,k}(\d s,\d z)= Z_{j,k}(s)\d sm(\d z) + |g_k(s,z) - g_j(s,z)|\d s\nu(\d z).
 \eeqnn
Let $\tilde{N}_{j,k}(\d s,\d z)= N_{j,k}(\d s,\d z) - \hat{N}_{j,k}(\d s,\d z)$ be the compensated random measure. We have
 \beqlb\label{s4.15}
Z_{j,k}(t) \ar=\ar M^c_{j,k}(t) + M^d_{j,k}(t) + \int_0^t [|\rho_j(s)-\rho_k(s)| - bZ_{j,k}(s)]\d s \cr
 \ar\ar\qquad
+ \int_0^t \d s \int_0^\infty |g_j(s,z) - g_k(s,z)|z \nu(\d z),
 \eeqlb
where $\{M^c_{j,k}(t): t\ge 0\}$ is a continuous local martingale with quadratic variation $2 cZ_{j,k}(t)\d t$ and
 \beqnn
M^d_{j,k}(t) = \int_0^t\int_0^\infty z \tilde{N}_{j,k}(\d s,\d z)
 \eeqnn
is a purely discontinuous local martingale. Using H\"{o}lder's inequality and Doob's martingale inequality we get
 \beqnn
\mbf{P}\Big[\sup_{0\le s\le t}Z_{j,k}(s)\Big]
 \ar\le\ar
\mbf{P}\bigg[\int_0^t |\rho_j(s)-\rho_k(s)|\d s\bigg] + |b|\mbf{P}\bigg[\int_0^t Z_{j,k}(s)\d s\bigg] \cr
 \ar\ar
+\,\mbf{P}\bigg[\int_0^t\d s\int_0^\infty |g_j(s,z) - g_k(s,z)|z \nu(\d z)\d s\bigg] \cr
 \ar\ar
+\, \mbf{P}\Big[\sup_{0\le s\le t} |M_{j,k}^c(t)|\Big] + \mbf{P}\bigg[\sup_{0\le s\le t}\bigg|\int_0^s\int_0^\infty z\tilde{N}(\d r,\d z)\bigg|\bigg] \cr
 \ar\le\ar
\|\rho_j-\rho_k\|_t+\|g_j-g_k\|_{\nu,t} + |b|\mbf{P}\bigg[\int_0^t Z_{j,k}(s)\d s\bigg]
\cr
 \ar\ar
+\, \mbf{P}\Big[\sup_{0\le s\le t} |M_{j,k}^c(t)|\Big] + \mbf{P}\bigg[\sup_{0\le s\le t}\bigg|\int_0^s\int_0^1 z\tilde{N}(\d r,\d z)\bigg|\bigg] \cr
 \ar\ar
+\, \mbf{P}\bigg[\bigg|\int_0^t\int_1^\infty z N(\d r,\d z)\bigg| + \bigg|\int_0^t\int_1^\infty z \hat{N}(\d r,\d z)\bigg|\bigg] \cr
 \ar\le\ar
\|\rho_j-\rho_k\|_t+\|g_j-g_k\|_{\nu,t} + |b|\mbf{P}\bigg[\int_0^t Z_{j,k}(s)\d s\bigg]
\cr
 \ar\ar
+\, 2\bigg\{c\mbf{P}\bigg[\int_0^t Z_{j,k}(s)\d s\bigg]\bigg\}^{\frac{1}{2}} + 2\bigg\{\mbf{P}\bigg[\int_0^t Z_{j,k}(s)\d s\int_0^1 u^2 m(\d u)\bigg]\bigg\}^{\frac{1}{2}} \cr
 \ar\ar
+\, 2\bigg\{\mbf{P}\bigg[\int_0^t\d s\int_0^1 |g_j(s,z) - g_k(s,z)|z^2 \nu(\d z)\bigg]\bigg\}^{\frac{1}{2}} \cr
 \ar\ar
+\,2\mbf{P}\bigg[\int_0^tZ_{j,k}(s)\d s \int_1^\infty um(\d u)\bigg] \cr
 \ar\ar
+\, 2\mbf{P}\bigg[\int_0^t\d s\int_1^\infty |g_j(s,z) - g_k(s,z)|z \nu(\d z)\bigg].
\eeqnn
By \eqref{s4.14} one can find a locally bounded function $t\mapsto C(t)$ so that
 \beqnn
\mbf{P}\Big[\sup_{0\le s\le t}Z_{j,k}(s)\Big]
 \le
C(t)\Big(\|\rho_j-\rho_k\|_t + \|g_j-g_k\|_{\nu,t} + \sqrt{\|\rho_j-\rho_k\|_t} + \sqrt{\|g_j-g_k\|_{\nu,t}}\Big).
 \eeqnn
It follows that
 \beqnn
\lim_{j,k\to\infty} \mbf{P}\Big[\sup_{0\le s\le t} |Y_j(s) - Y_k(s)|\Big]
 \le
\lim_{j,k\to\infty} \mbf{P}\Big[\sup_{0\le s\le t} Z_{j,k}(s)\Big] = 0.
 \eeqnn
Then there is an increasing sequence of integers $\{k_n\}$ such that
 \beqlb\label{s4.16}
\mbf{P}\Big[\sup_{0\le s\le n}|Y_j(s) - Y_k(s)|\Big]\le 1/2^n, \qquad j,k\ge k_n.
 \eeqlb
Consequently, for every $t\ge 0$,
 \beqnn
\mbf{P}\bigg[\sum_{n=1}^\infty\sup_{0\le s\le t}|Y_{k_n}(s) - Y_{k_{n+1}}(s)|\bigg]< \infty.
 \eeqnn
Thus we have a.s.\
 \beqnn
\sum_{n=1}^\infty\sup_{0\le s\le t}|Y_{k_n}(s) - Y_{k_{n+1}}(s)|< \infty.
 \eeqnn
That yields the existence of a c\`{a}dl\`{a}g positive process $\{Y(t): t\ge 0\}$ so that \eqref{s4.12} holds a.s.\ for every $t\ge 0$. By letting $j\to \infty$ along the sequence $\{k_n\}$ in \eqref{s4.16} we get
 \beqnn
\mbf{P}\Big[\sup_{0\le s\le n}|Y(s) - Y_k(s)|\Big]\le 1/2^n, \qquad k\ge k_n.
 \eeqnn
Then we get \eqref{s4.11}. \qed

\noindent{\textit{Proof of Theorem~\ref{ts4.3}.~~}} By Proposition~10.3 in Li (2011, p.236), there are positive sequences $\{\rho_k\}\subset \mcr{L}^0$ and $\{g_k\}\subset \mcr{L}^0_\nu(0,\infty)$ such that $d(\rho_k,\rho) + d_\nu(g_k,g)\to 0$ as $k\to \infty$. Let $\{Y_k(t): t\ge 0\}$ be defined by \eqref{s4.7} with $\rho= \rho_k$ and $g= g_k$. By Lemma~\ref{ts4.5}, there is a positive c\`{a}dl\`{a}g process $\{Y(t): t\ge 0\}$ and a sequence $\{k_n\}$ so that \eqref{s4.11} and \eqref{s4.12} hold for every $t\ge 0$. For each $k\ge 1$ we have $|Y_t-Y_k(t)|\le Z_k(t)$, where
 \beqnn
Z_k(t) \ar=\ar \int_0^t h_{t-s}|\rho(s)-\rho_k(s)|\d s + \int_0^t \int_{\rho(s)\land \rho_k(s)}^{\rho(s)\vee \rho_k(s)} \int_W w(t-s) N_0(\d s,\d u,\d w) \cr
 \ar\ar\qquad
+ \int_0^t\int_0^\infty\int_{g(s,z)\land g_k(s,z)}^{g(s,z)\vee g_k(s,z)} \int_W w(t-s) N_1(\d s,\d z,\d u,\d w).
 \eeqnn
By the arguments leading to \eqref{s4.14} we get
 \beqnn
\mbf{P}[|Y_t-Y_k(t)|]\le \mbf{P}[Z_k(t)]
 \le
\e^{|b|t}(\|\rho-\rho_k\|_t+\|g-g_k\|_{\nu,t}).
 \eeqnn
It follows that $\mbf{P}[|Y_t-Y_k(t)|]\to 0$ as $n\to \infty$. By choosing a smaller sequence $\{k_n\}$ we have a.s.\ $Y_t= \lim_{n\to \infty} Y_{k_n}(t)= Y(t)$. Then the positive c\`{a}dl\`{a}g process $\{Y(t): t\ge 0\}$ is a modification of $\{Y_t: t\ge 0\}$. By Lemma~\ref{ts4.4}, property (4.B) holds for $\{Y_k(t): t\ge 0\}$ with $\rho=\rho_k$ and $g=g_k$. Then the property holds for $\{Y(t): t\ge 0\}$. By Proposition~\ref{ts4.1} the process $\{Y(t): t\ge 0\}$ also has property (4.A). \qed

\bgproposition\label{ts4.6} Suppose that $\rho_k,\rho\in \mcr{L}^1$ and $g_k, g\in \mcr{L}^1_\nu(0,\infty)$ are positive processes such that $d(\rho_k,\rho) + d_\nu(g_k,g)\to 0$ as $k\to \infty$. Let $\{Y_t: t\ge 0\}$ be the positive c\`{a}dl\`{a}g process defined by (\ref{s4.7}). Let $\{Y_k(t): t\ge 0\}$ be the positive c\`{a}dl\`{a}g process defined by the same formula with $\rho=\rho_k$. Then we have
 \beqlb\label{s4.17}
\lim_{k\to\infty} \mbf{P}\Big[\sup_{0\le s\le t} |Y_k(s)-Y_s|\Big] = 0, \qquad t\ge 0,
 \eeqlb
and there is a subsequence $\{k_n\}\subset \{k\}$ so that a.s.\
 \beqlb\label{s4.18}
\lim_{n\to \infty} \sup_{0\le s\le t}|Y_{k_n}(s)-Y_s| = 0, \qquad t\ge 0.
 \eeqlb
\edproposition

\proof For $j,k\ge 1$ let $\{Z_{j,k}(t): t\ge 0\}$ be defined as in \eqref{s4.13}. Then $\{(Z_{j,k}(t), \mcr{G}_t): t\ge 0\}$ is a CBI-process with predictable immigration rates given by $\{(|\rho_j(s) - \rho_k(s)|, |g_j(s,z) - g_k(s,z)|): s\ge 0, z>0\}$. By Proposition~\ref{ts4.1} and Theorem~\ref{ts4.3}, the properties (4.A) and (4.B) hold for $\{Z_{j,k}(t): t\ge 0\}$. In particular, a decomposition like \eqref{s4.15} is valid for this process. The remaining arguments go as in the proof of Lemma~\ref{ts4.5}. \qed

\section{Solutions of the stochastic equations}

 \setcounter{equation}{0}

In this section, we prove there is a pathwise unique solution to \eqref{s1.13} and show the solution is a Markov process with generator $(L,\mcr{D})$. Suppose that $\{(X_t,\mcr{F}_t)\}$, $\{N_0(\d s,\d u,\d w)\}$ and $\{N_1(\d s,\d z,\d u,\d w)\}$ are given as in the introduction. Let the filtration $(\mcr{G}_t)$ be defined as in Theorem~\ref{ts3.2}. The following result gives a reformulation of \eqref{s1.10} in terms of a martingale problem.

\bgproposition\label{ts5.1} A positive c\`{a}dl\`{a}g process $\{Y_t: t\ge 0\}$ is a weak solution of \eqref{s1.10} if and only if it solves the martingale problem of $(L,\mcr{D})$, that is, for every $f\in \mcr{D}$,
 \beqlb\label{s5.1}
f(Y_t) = f(Y_0) + \int_0^t Lf(Y_s)\d s + \mbox{local mart.}
 \eeqlb
\edproposition

\proof If $\{Y_t: t\ge 0\}$ is a weak solution of \eqref{s1.10}, we may use It\^{o}'s formula to see it solves the martingale problem given by \eqref{s5.1}. Conversely, let us assume $\{Y_t: t\ge 0\}$ is a solution of the martingale problem \eqref{s5.1}. By Proposition~\ref{ts4.1}, the process has property (4.A) in Section~4 with $\rho(s) = \beta(Y_{s-})$ and $g(s,z) = q(Y_{s-},z)$. By Theorem~III.7.1$'$ in Ikeda and Watanabe (1989, p.90), on an extension of the original probability space there is a Brownian motion $\{B(s)\}$ so that
 \beqnn
M^c(t) =\int_0^t\sqrt{2cY_{s-}} \d B(s), \qquad t\ge 0.
 \eeqnn
By Theorem~III.7.4 in Ikeda and Watanabe (1989, p.93), on a further extension of the probability space there are independent Poisson time-space random measures $M(\d s,\d z,\d u)$ and $N(\d s,\d z)$ with intensities $\d sm(\d z)\d u$ and $\d s\nu(\d z)$, respectively, so that
 \beqnn
\int_0^t \int_0^\infty z \tilde{N}_0(\d s,\d z)
 =
\int_0^t\int_0^\infty \int_0^{Y_{s-}} z \tilde{M}(\d s,\d z,\d u) +
\int_0^t\int_0^\infty\int_0^{g(Y_{s-},z)} z N(\d s,\d z,\d u).
 \eeqnn
Then $\{Y_t\}$ is a weak solution of the stochastic equation \eqref{s1.10}. \qed

\bgproposition\label{ts5.2} Let $\{Y_t: t\ge 0\}$ be a solution to \eqref{s1.13}. Then there is a locally bounded function $t\mapsto C(t)$ so that
 \beqlb\label{s5.2}
\mbf{P}[Y_t]
 \le
C(t)\big(1+ \mbf{P}[X_0]+ \sqrt{\mbf{P}[X_0]}\,\big), \qquad t\ge 0.
 \eeqlb
\edproposition

\proof Let $\{\tau_k: k\ge 1\}$ be the increasing sequence of stopping times defined by $\tau_k= \inf\{t\ge 0: Y_t\ge k\}$. Then we have a.s.\ $\lim_{k\to \infty}\tau_k= \infty$. By \eqref{s1.13} and Condition~(1.A) we have
 \beqnn
\mbf{P}[Y_t1_{\{t<\tau_k\}}] \ar\le\ar \mbf{P}[X_{t\land\tau_k}] + \mbf{P}\bigg[\int_0^{t\land\tau_k}\big(h_{t-s} + \mbf{N}_0[w(t-s)]\big)\beta(Y_{s-})\d s\bigg] \cr
 \ar\ar\qquad
+\, \mbf{P}\bigg[\int_0^{t\land\tau_k}\d s\int_0^\infty\mbf{Q}_z[w(t-s)]q(Y_{s-},z) \nu(\d z)\bigg] \cr
 \ar\le\ar
\mbf{P}[X_{t\land\tau_k}] + \mbf{P}\bigg[\int_0^{t\land\tau_k} \e^{-b(t-s)}\bigg(\beta(Y_{s-}) + \int_0^\infty q(Y_{s-},z)z \nu(\d z)\bigg)\d s\bigg] \cr
 \ar\le\ar
\mbf{P}\Big[\sup_{0\le s\le t}X_s\Big] + K\mbf{P}\bigg[\int_0^t \e^{-b(t-s)}(1+Y_{s-})1_{\{s<\tau_k\}}\d s\bigg].
 \eeqnn
By Corollary~7.3 in Li (2018b), there is a locally bounded function $t\mapsto C_0(t)$ so that
 \beqnn
\mbf{P}\Big[\sup_{0\le s\le t}X_s\Big]\le C_0(t)\big(\mbf{P}[X_0] + \sqrt{\mbf{P}[X_0]}\,\big).
 \eeqnn
Then $t\mapsto \mbf{P}[Y_t1_{\{t<\tau_k\}}]$ is locally bounded and
 \beqnn
\mbf{P}[Y_t1_{\{t<\tau_k\}}] \ar=\ar  C_0(t)\big(\mbf{P}[X_0] + \sqrt{\mbf{P}[X_0]}\,\big) + K\mbf{P}\bigg[\int_0^t \e^{-b(t-s)}(1+Y_s) 1_{\{s< \tau_k\}}\d s\bigg] \cr
 \ar\le\ar
C_0(t)\big(\mbf{P}[X_0] + \sqrt{\mbf{P}[X_0]}\,\big) + Kt\e^{|b|t} + K\e^{|b|t}\int_0^t \mbf{P}[Y_s1_{\{s< \tau_k\}}]\d s.
 \eeqnn
By Gronwall's inequality, we can can find a function $t\mapsto C(t)$ independent of $k\ge 1$ so that
 \beqlb\label{s5.3}
\mbf{P}[Y_t1_{\{t< \tau_k\}}]
 \le
C(t)\big(1+ \mbf{P}[X_0]+ \sqrt{\mbf{P}[X_0]}\,\big).
 \eeqlb
Then \eqref{s5.2} follows from \eqref{s5.3} by Fatou's lemma.\qed

\bgproposition\label{ts5.3} There is at most one solution to \eqref{s1.13}. \edproposition

\proof Suppose that $\{Y_t: t\ge 0\}$ and $\{Z_t: t\ge 0\}$ are two solutions of the equation. Then we have $|Y_t - Z_t|\le \xi_t\le Y_t + Z_t$, where
 \beqlb\label{s5.4}
\xi_t \ar=\ar \int_0^t h_{t-s}|\beta(Y_{s-})-\beta(Z_{s-})|\d s + \int_0^t \int_{\beta(Y_{s-})\land \beta(Z_{s-})}^{\beta(Y_{s-})\vee \beta(Z_{s-})} \int_W w(t-s) N_0(\d s,\d u,\d w) \cr
 \ar\ar\qqquad
+ \int_0^t \int_0^\infty \int_{q(Y_{s-},z)\land q(Z_{s-},z)}^{q(Y_{s-},z)\vee q(Z_{s-},z)}\int_W w(t-s) N_1(\d s,\d z,\d u,\d w).
 \eeqlb
By Proposition~\ref{ts5.2}, the function $t\mapsto \mbf{P}[Y_t + Z_t]= \mbf{P}[Y_t] + \mbf{P}[Z_t]$ is locally bounded, then so is $t\mapsto \mbf{P}[\xi_t]$. We can rewrite \eqref{s5.4} into
 \beqnn
\xi_t \ar=\ar \int_0^t h_{t-s}|\beta(Y_{s-})-\beta(Z_{s-})|\d s + \int_0^t \int_0^{|\beta(Y_{s-})-\beta(Z_{s-})|} \int_W w(t-s) M_0(\d s,\d u,\d w) \cr
 \ar\ar\qqquad
+ \int_0^t \int_0^\infty \int_0^{|q(Y_{s-},z)-q(Z_{s-},z)|}\int_W w(t-s) M_1(\d s,\d z,\d u,\d w),
 \eeqnn
where $M_0(\d s,\d u,\d w)$ is a spatial shift of $N_0(\d s,\d u,\d w)$ and $M_0(\d s,\d u,\d w)$ is a spatial shift of $N_1(\d s,\d z,\d u,\d w)$. By Proposition~\ref{ts4.2} and Condition~(1.B),
 \beqlb\label{s5.5}
\mbf{P}[\xi_t]
 \ar=\ar
\mbf{P}\bigg[\int_0^t \e^{-b(t-s)}\bigg(|\beta(Y_{s-})-\beta(Z_{s-})| + \int_0^\infty |q(Y_{s-},z)-q(Z_{s-},z)|z \nu(\d z)\bigg)\d s\bigg] \cr
 \ar\le\ar
\e^{|b|t}\mbf{P}\bigg[\int_0^t r(|Y_{s-}-Z_{s-}|) \d s\bigg]
 \le
\e^{|b|t}\int_0^t r(\mbf{P}[|Y_s-Z_s|]) \d s \cr
 \ar\le\ar
\e^{|b|t}\int_0^t r(\mbf{P}[\xi_s]) \d s,
 \eeqlb
where the last two inequalities hold by the convexity and monotonicity of $u\mapsto r(u)$. Now define the continuous increasing function
 \beqnn
u(t)= \int_0^t r(\mbf{P}[\xi_s]) \d s, \qquad t\ge 0.
 \eeqnn
We claim that $u(t)= 0$ for every $t\ge 0$. Otherwise, there is $a>0$ so that $u(a)>0$. Let $v= \inf\{t\in [0,a]: u(t)>0\}$. Then $0\le v<a$ and $u(v)=0$. Since $u\mapsto r(u)$ is increasing, we get $\d u(s)= r(\mbf{P}[\xi_s])\d s\le r(\e^{|b|s}u(s))\d s$ by \eqref{s5.5}. For $v<t \le a$ we have
 \beqnn
a-t\ge \int_t^a{r(\e^{|b|s}u(s))\d s\over r(\e^{|b|a}u(s))}
 \ge
\int_t^a {\d u(s)\over r(\e^{|b|a}u(s))}
 =
\e^{-|b|a}\int_{\e^{|b|a}u(t)}^{\e^{|b|a}u(a)} {\d u\over r(u)}.
 \eeqnn
By letting $t\to v$ and using Condition~(1.B) we conclude
 \beqnn
a-v\ge \e^{-|b|a}\int_0^{\e^{|b|a} u(a)} r(u)^{-1}\d u= \infty,
 \eeqnn
which gives a contradiction. Then we must have $u(t)=0$ and hence $\mbf{P}[\xi_t]= 0$ for every $t\ge 0$. That proves the pathwise uniqueness of solution to \eqref{s1.13}. \qed

\bgtheorem\label{ts5.4} There is a pathwise unique solution $\{Y_t: t\ge 0\}$ to \eqref{s1.13}. Moreover, the process $\{(Y_t,\mcr{G}_t): t\ge 0\}$ solves the martingale problem \eqref{s5.1} for $(L,\mcr{D})$.
\edtheorem

\proof The pathwise uniqueness of solution to \eqref{s1.13} follows from Proposition~\ref{ts5.3}. We shall give the construction of a solution to the stochastic equation in two steps.

\textit{Step~1.} Instead of {\rm(1.A)}, we assume the following stronger condition: there is a constant $K\ge 0$ so that
 \beqlb\label{s5.6}
\beta(x) + \int_{(0,\infty)} q(x,z)z\nu(\d z)\le K, \qquad x\ge 0.
 \eeqlb
Let $Y_0(t) = X_t$ and define inductively
 \beqnn
\rho_k(s)= \beta(Y_{k-1}(s-)), ~ g_k(s,z)= q(Y_{k-1}(s-),z)
 \eeqnn
and
 \beqnn
Y_k(t) \ar=\ar X_t + \int_0^t h_{t-s}\rho_k(s)\d s + \int_0^t \int_0^{\rho_k(s)} \int_W w(t-s) N_0(\d s,\d u,\d w) \cr
 \ar\ar\qquad
+ \int_0^t\int_0^\infty\int_0^{g_k(s,z)} \int_W w(t-s) N_1(\d s,\d z,\d u,\d w).
 \eeqnn
Then by Proposition~\ref{ts4.2} we one can show
 \beqlb\label{s5.7}
\mbf{P}[Y_k(t)]\le C_1(t):= \e^{-bt}\mbf{P}[X_0] + K\int_0^t \e^{-b(t-s)}\d s.
 \eeqlb
For $k,j\ge 1$ let the process $\{Z_{j,k}(t): t\ge 0\}$ be defined as in \eqref{s4.13}. Then $|Y_j(t) - Y_k(t)|\le Z_{j,k}(t)\le Y_j(t) + Y_k(t)$. From \eqref{s5.7} it follows that $\mbf{P}[Z_{j,k}(t)]\le 2C_1(t)$ for $t\ge 0$. As in \eqref{s5.5} one can see
 \beqlb\label{s5.8}
\mbf{P}[Z_{j,k}(t)]
 \le
\e^{|b|t}\int_0^t r\big(\mbf{P}[|Y_{j-1}(s)-Y_{k-1}(s)|]\big) \d s
 \le
\e^{|b|t}\int_0^t r\big(\mbf{P}[Z_{j-1,k-1}(s)]\big) \d s.
 \eeqlb
Let $R_n(t)= \sup_{j,k\ge n}\mbf{P}[Z_{j,k}(t)]\le 2C_1(t)$. By \eqref{s5.8} and dominated convergence,
 \beqnn
\lim_{n\to \infty} R_n(t)
 \le
\e^{|b|t}\int_0^t r\Big(\lim_{n\to \infty} R_{n-1}(s)\Big) \d s
 =
\e^{|b|t}\int_0^t r\Big(\lim_{n\to \infty} R_n(s)\Big) \d s.
 \eeqnn
As in the last part of the proof of Proposition~\ref{ts5.3} we see
 \beqlb\label{s5.9}
\lim_{j,k\to \infty}\mbf{P}[Z_{j,k}(t)] =\lim_{n\to \infty} R_n(t) = 0.
 \eeqlb
Observe that
 \beqnn
\|\rho_j-\rho_k\|_t +\|g_j-g_k\|_{t,\nu}
 \ar=\ar
\mbf{P}\bigg[\int_0^t\bigg(|\beta(Y_{j-1}(s-))-\beta(Y_{k-1}(s-))| \cr
 \ar\ar
+ \int_0^\infty |q(Y_{j-1}(s-),z)-q(Y_{k-1}(s-),z)|z\nu(\d z)\bigg)\d s\bigg] \cr
 \ar\le\ar
\e^{|b|t}\mbf{P}\bigg[\int_0^t r\big(|Y_{j-1}(s)-Y_{k-1}(s)|\big)\d s\bigg] \cr
 \ar\le\ar
\e^{|b|t}\int_0^t r\big(\mbf{P}[Z_{j-1,k-1}(s)]\big)\d s.
 \eeqnn
Then \eqref{s5.9} implies that $\|\rho_j-\rho_k\|_t +\|g_j-g_k\|_{t,\nu}\to 0$ as $j,k\to \infty$. Since $\mcr{L}^1$ and $\mcr{L}^1_\nu(0,\infty)$ are complete, there are positive processes $\rho\in \mcr{L}^1$ and $g\in \mcr{L}^1_\nu(0,\infty)$ so that $\|\rho_k-\rho\|_t +\|g_k-g\|_{t,\nu}\to 0$ as $k\to \infty$. Define the process $\{Y_t: t\ge 0\}$ by \eqref{s4.7}. By Condition~(1.B) we have
 \beqnn
\ar\ar\int_0^t\mbf{P}\bigg[|\rho_k(s)-\beta(Y_{s-})| + \int_{(0,\infty)} |g_k(s,z)-q(Y_{s-},z)|z \nu(\d z)\bigg]\d s \cr
 \ar\ar\qquad
= \int_0^t\mbf{P}\bigg[|\beta(Y_{k-1}(s))-\beta(Y_s)| + \int_{(0,\infty)} |q(Y_{k-1}(s),z) - q(Y_s,z)|z \nu(\d z)\bigg]\d s \cr
 \ar\ar\qquad
\le \int_0^t\mbf{P}\big[r(|Y_{k-1}(s)-Y_s|)\big]\d s
 \le
\int_0^t r\big(\mbf{P}[|Y_{k-1}(s)-Y_s|]\big)\d s.
 \eeqnn
By Proposition~\ref{ts4.6}, the right-hand side vanishes as $k\to \infty$. Then we may identify $s\mapsto \rho(s)$ and $s\mapsto \beta(Y(s-))$ as elements of $\mcr{L}^1$ and identify $(s,z)\mapsto g(s,z)$ and $(s,z)\mapsto q(Y(s-),z)$ as elements of $\mcr{L}^1_\nu(0,\infty)$. Now \eqref{s1.13} follows from \eqref{s4.7}. The martingale problem characterization \eqref{s5.1} follows by Theorem~\ref{ts4.3}.

\textit{Step~2.} In the general case where \eqref{s5.6} is not necessarily true, for $n\ge 1$ we consider the stochastic equation
 \beqlb\label{s5.10}
Y_t \ar=\ar X_t + \int_0^t h_{t-s}\beta(Y_{s-})\d s + \int_0^t\int_0^{\beta(Y_{s-}\land n)} \int_W w(t-s) N_0(\d s,\d u,\d w) \cr
 \ar\ar\qquad
+ \int_0^t\int_0^\infty\int_0^{q(Y_{s-}\land n,z)} \int_W w(t-s) N_1(\d s,\d z,\d u,\d w).
 \eeqlb
By Step~1 and Proposition~\ref{ts5.3}, there is a pathwise unique solution $\{Y_n(t): t\ge 0\}$ to \eqref{s5.10}. Let $\tau_n= \inf\{t\ge 0: Y_n(t)\ge n\}$. Then $Y_n(t)= Y_{n+1}(t)$ for $0\le t<\tau_n$. As in the proof of Proposition~\ref{ts5.2}, one can see $\tau_n\to \infty$ increasingly in probability as $t\to \infty$. Then $Y(t)=\lim_{n\to \infty}Y_k(t)$ defines a process $\{Y(t): t\ge 0\}$. From \eqref{s5.10} we see this process is a solution to \eqref{s1.13}. \qed

Comparison properties of stochastic equations of the form \eqref{s1.10} were studied in Bertoin and Le~Gall (2006), Dawson and Li (2012) and Fu and Li (2010). Those results have played important roles in the study of stochastic flows induced by those equations. The next theorem provides a comparison result for the stochastic equation \eqref{s1.13}.

\bgtheorem\label{ts5.5} Let $(\beta,q)$ and $(\beta^\prime,q^\prime)$ be two sets of parameters satisfying Conditions~{\rm(1.A)} and~{\rm(1.B)}. Suppose that:
 \bitemize

\itm $x\mapsto \beta(x)$ or $x\mapsto \beta^\prime(x)$ is increasing on $[0,\infty)$;

\itm for each $z>0$, $x\mapsto q(x,z)$ or $x\mapsto q^\prime(x,z)$ is increasing on $[0,\infty)$.

\itm $\beta(x)\le \beta^\prime(x)$ and $q(x,z)\le q^\prime(x,z)$ for all $x\ge 0$ and $z>0$.

 \eitemize
Let $\{Y_t: t\ge 0\}$ be the solution of \eqref{s1.13} and $\{Y_t^\prime: t\ge 0\}$ the solution of the stochastic equation with $(\beta,q)$ replaced by $(\beta^\prime,q^\prime)$. Then $\mbf{P}(Y_t\le Y_t^\prime$ for all $t\ge 0) = 1$.
\edtheorem

\proof For simplicity, we assume $x\mapsto \beta(x)$ is increasing. Let $A=\{z>0: x\mapsto q(x,z)$ is increasing$\}$. Then $x\mapsto q^\prime(x,z)$ is increasing for $z\in A^\prime:= (0,\infty)\setminus A$. Let $\{Y_k(t)\}$ be the sequence defined as in the proof of Theorem~\ref{ts5.4}, and let $\{Y_k^\prime(t)\}$ be defined in the same way with $(\beta,q)$ replaced by $(\beta^\prime,q^\prime)$. By induction in $k\ge 1$ we see that
 \beqnn
Y_k(t) \ar=\ar X_t + \int_0^t h_{t-s}\beta(Y_{k-1}(s-))\d s + \int_0^t \int_0^{\beta(Y_{k-1}(s-))} \int_W w(t-s) N_0(\d s,\d u,\d w) \cr
 \ar\ar\qqquad
+ \int_0^t\int_A\int_0^{q(Y_{k-1}(s-),z)} \int_W w(t-s) N_1(\d s,\d z,\d u,\d w) \cr
 \ar\ar\qqquad
+ \int_0^t\int_{A^\prime}\int_0^{q^\prime(Y_{k-1}(s-),z)} \int_W w(t-s) N_1(\d s,\d z,\d u,\d w) \cr
 \ar\le\ar
X_t + \int_0^t h_{t-s}\beta(Y_{k-1}^\prime(s-))\d s + \int_0^t \int_0^{\beta(Y_{k-1}^\prime(s-))} \int_W w(t-s) N_0(\d s,\d u,\d w) \cr
 \ar\ar\qqquad
+ \int_0^t\int_A\int_0^{q(Y_{k-1}^\prime(s-),z)} \int_W w(t-s) N_1(\d s,\d z,\d u,\d w) \cr
 \ar\ar\qqquad
+ \int_0^t\int_{A^\prime}\int_0^{q^\prime(Y_{k-1}^\prime(s-),z)} \int_W w(t-s) N_1(\d s,\d z,\d u,\d w) \cr
 \ar\le\ar
X_t + \int_0^t h_{t-s}\beta^\prime(Y_{k-1}^\prime(s-))\d s + \int_0^t \int_0^{\beta^\prime(Y_{k-1}^\prime(s-))} \int_W w(t-s) N_0(\d s,\d u,\d w) \cr
 \ar\ar\qquad
+ \int_0^t\int_0^\infty\int_0^{q^\prime(Y_{k-1}^\prime(s-),z)} \int_W w(t-s) N_1(\d s,\d z,\d u,\d w) \cr
 \ar=\ar
Y_k^\prime(t).
 \eeqnn
Then $\mbf{P}(Y_t\le Y_t^\prime$ for all $t\ge 0) = 1$ by the proof of Theorem~\ref{ts5.4}. \qed

\bigskip

\textbf{Acknowledgments} ~ I would like to thank Rongjuan Fang and an anonymous referee for their careful reading of the paper and very helpful comments on the presentation of the results. I am grateful to the Laboratory of Mathematics and Complex Systems (Ministry of Education) for providing the research facilities to carry out the project.

 \end{document}